\documentclass[11pt]{article}
\usepackage{mathrsfs}

\usepackage{amsmath}
\usepackage{amssymb}
\usepackage{amsfonts}
\usepackage{amsthm}

\usepackage{latexsym}
\usepackage{graphicx}

\renewcommand{\epsilon}{\varepsilon}

\newtheorem{example}{Example}
\newtheorem{remark}{Remark}

\def\expect{{\mathbb  E}}
\def\var{{\mathbb Var}}
\def\Pr{{\mathbb P}}

\def\eqdef{\triangleq}
\def\real{{\mathbb  R}}
\def\nat{{\mathbb  N}}
\def\eqd{\stackrel{d}{=}}
\def\ind{{\bf 1}}
\def\leqd{\stackrel{d}{\leq}}
\def\12{\frac{1}{2}}

\def\mp{{\Pi}} 
\def\rmp{{M}} 
\def\mdp{{J}}
\def\rbp{{\Lambda}} 

\newtheorem{theorem}{Theorem}

\newtheorem{lemma}{Lemma}

\newtheorem{corollary}{Corollary}
\newtheorem{definition}{Definition}

\begin{document}

\author{Predrag R. Jelenkovi\'c \hspace{0.3in} Jian Tan  \\
\begin{tabular}{ c }
\small Department of Electrical Engineering \\
\small Columbia University, New York, NY 10027 \\
\small \{predrag, jiantan\}@ee.columbia.edu
\end{tabular}
}

\title{ \textbf{Modulated Branching Processes, Origins of Power Laws
and Queueing Duality}}

\date{\small September 2007; revised July 2008}

\maketitle

\thispagestyle{empty}

\footnotetext{The preliminary version of this work has appeared in
the extended abstract in Proceedings of the  }
\footnotetext{Forty-Fourth Annual
Allerton Conference, 
Allerton House, UIUC, Illinois, USA, September 2006.}
\footnotetext{Technical Report EE2007-09-25, Department of
Electrical Engineering, Columbia University,}
 \footnotetext{
  New York, NY,  September 25, 2007.}
\begin{abstract}

Power law distributions have been repeatedly observed in a wide
variety of socioeconomic, biological and technological areas. In
many of the observations, e.g., city populations and sizes of living
organisms, the objects of interest evolve due to the replication of
their many independent components, e.g., births-deaths of
individuals and replications of cells. Furthermore, the rates of the
replication are often controlled by exogenous parameters causing
periods of expansion and contraction, e.g., baby booms and busts,
economic booms and recessions, etc. In addition, the sizes of these
objects often have reflective lower boundaries, e.g., cities do not
fall bellow a certain size, low income individuals are subsidized by
the government, companies are protected by bankruptcy laws, etc.

Hence, it is natural to propose reflected modulated branching
processes as generic models for many of the preceding observations.
Indeed, our main results show that the proposed mathematical models
result in power law distributions under quite general {polynomial
G\"artner-Ellis} conditions,  the generality of which could explain
the ubiquitous nature of power law distributions. In addition, on a
logarithmic scale, we establish an asymptotic equivalence between
the reflected branching processes and the corresponding
multiplicative ones. The latter, as recognized by Goldie (1991)
\cite{Goldie91}, is known to be dual to queueing/additive processes.
We emphasize this duality further in the generality of stationary
and ergodic processes.

 \vspace{0.3in}
 \noindent {\bf Keywords:} Modulated branching processes, reflective/absorbing barriers,
 reflected multiplicative processes,
proportional growth models, power law distributions, heavy tails,
subexponential distributions, queueing processes,  reflected
additive random walks, Cram\'er large deviations, polynomial
G\"artner-Ellis conditions.

\end{abstract}

\newpage

\section{Introduction}

Power law distributions are found in a wide range of domains,
ranging from socioeconomic to biological and technological areas.
Specifically, these types of distributions describe the city
populations, species-area relationships, sizes of living organisms,
values of companies, distributions of wealth, and more recently,
sizes of documents on the Web,  visitor access patterns on Web
sites, etc. Hence, one would expect that there exist universal
mathematical laws that explain this ubiquitous nature of power law
distributions. To this end, we propose a class of models, termed
modulated branching processes with reflective lower barriers that,
under quite general \emph{polynomial G\"artner-Ellis} conditions,
result in power law distributions.

Empirical observations of power laws have a long history,  starting
from the discovery by Pareto \cite{Vi97} in 1897 that a plot of the
logarithm of the
  number of incomes above a level against the logarithm of that level
  yields points close to a straight line, which is essentially equivalent
  to saying that the income distribution follows a power law.
Hence, power law distributions are often called Pareto
distributions; for more recent study on income distributions see
\cite{Ch53,B60,C96,R01,R03}. In a different context, early work by
Arrhenius \cite{O21} in 1921 conjectured a power law relationship
between the number of species and the census area, which was
followed by Preston's prediction in \cite{Pr62} that the slope on
the log/log species-area plot has a canonical value equal to
$0.262$; for additional information and measurements on species-area
relationships see \cite{CM79,Pl00,Ke03}. Interestingly, there also
exists a power law relationship between the rank of the cities and
the population of the corresponding cities. This was proposed by
Auerbach \cite{Au13} in 1913 and later studied by Zipf \cite{Zi49}
in 1949, after whom power law is also known as Zipf's law. Ever
since, much attention on both empirical examinations and
explanations of the city size distributions have been drawn
\cite{Zi49,IO03,Ga99,RR80,Pa85,ALC01}. Similar observations have
been made for firm sizes \cite{ABH97}, language family sizes
\cite{wich05}, and even the gene family and protein statistics
\cite{HN98,RzG01,LQZ02,BHWF06}.
It is maybe even more surprising that many features of the Internet
 are governed by power laws, including the distribution of pages per
 Web site \cite{AH99}, the page request distribution
 \cite{CBC95,BC99}, the file size distribution \cite{Do01,JelenkovicHetero03},
 Ethernet LAN traffic \cite{Leland93}, World Wide Web traffic
 \cite{CB97}, the number
 of visitors per Web site \cite{HA00,AH02}, the distribution of scenes in MPEG video
 streams \cite{JLS97}
  and the distribution of the
  indegrees and outdegrees in the Web graph as well as the physical network connectivity graph
  \cite{FFF99, BAJ99, KRRST00, MMB00}.
In socio-economic areas, in addition to income distributions,  the
fluctuations in stock prices have also been observed to be
characterized by power laws \cite{GGP03, LS96a2}. This paragraph
only exemplifies various observations of power laws; for a more
complete survey see \cite{Mi04}.

Hence, these repeated empirical observations of power laws over a
period of more than a hundred years strongly suggest that there
exist general mathematical laws that govern these phenomena. In this
regard, after carefully examining the situations that result in
power laws, we discover that most of them are characterized by the
following three features. First, in the vast majority of these
observations, e.g., city populations and sizes of living organisms,
the objects of interest evolve due to the replication of their many
independent components, e.g., birth-deaths of individuals and
replications of cells. Secondly, the rate of replication of the many
components is often controlled by exogenous parameters causing
periods of baby booms and busts, economic growths and recessions,
etc. Thirdly, the sizes of these objects often have lower
boundaries, e.g., cities do not fall bellow a certain size, low
income individuals are subsidized by the government, companies are
protected by bankruptcy laws, etc.

In order to capture the preceding features, it is natural to propose
{\it modulated branching processes} (MBP) with reflective or
absorbing barriers as generic models for many of the observations of
power laws. Indeed, one of our main results, presented in
Theorem~\ref{theorem:rbp}, shows that MBPs with reflective barriers
almost invariably produce power law distributions under quite
general \emph{polynomial G\"artner-Ellis} conditions. The generality
of our results could explain the ubiquitous nature of power law
distributions. Furthermore, an informal interpretation of our main
results, stated in Theorems~\ref{theorem:rbp} and
\ref{theorem:lessOne} of Section~\ref{s:main}, suggests that
alternating periods of expansions and contractions, e.g., economic
booms and recessions, are primarily responsible for the appearance
of power law distributions. Actually, Theorem \ref{theorem:lessOne}
shows that the distribution of the reflected MBP is exponentially
bounded if the process always contracts. From a mathematical
perspective, we develop a novel sample path technique for analyzing
reflected modulated branching processes since these objects appear
new and the traditional methods for investigating branching
processes \cite{AV72} do not directly apply; a preliminary version
of this work has appeared in the extended abstract in
\cite{PJ06MBPconf}.

Formal description of our reflected modulated branching process
(RMBP) model is given in Section~\ref{s:model}. In the singular case
when the number of individuals born in each state of the modulating
process is constant, our model reduces to a reflected multiplicative
process. A rigorous connection (duality) between the reflected
multiplicative processes (RMPs) and queueing theory was established
in Section 5 of Goldie (1991) \cite{Goldie91}; this duality was
repeatedly observed and used later in, e.g., \cite{SC97,GL05}. In
Subsection~\ref{ss:multiplicative} we further emphasize this duality
in the context of stationary and ergodic processes. We would like to
point out that this duality makes a vast literature on queueing
theory directly applicable to the analysis of RMPs. As a direct
consequence of this connection, in
Subsection~\ref{ss:multiplicative} we translate several well known
queueing results to the context of RMPs. Informally, these results
show that the role which exponential distributions play in queueing
theory, and in additive reflected random walks in general, is
represented by power low distributions in the framework of
RMPs/RMBPs. Furthermore, this relationship appears to reduce the
debate on the relative importance of power law versus exponential
distributions/models to the analogous question of the prevalence of
proportional growth versus additive phenomena. Interestingly, the
power law distribution satisfies the \emph{memoryless property} in
the multiplicative world, playing an equivalent role to the
memoryless exponential distribution in the additive world. Indeed,
if $\Pr[M>x]=x^{-\alpha}, \alpha>0, x\geq 1,$ then, for $x, y\geq
1$, we obtain $\Pr[M>xy | M>x]= \Pr[M>y]$.

Furthermore, this duality immediately implies and generalizes many
of the prior results in the area of RMPs and power laws, e.g., see
\cite{LS96b,LS96a1,LS96a2,SC97}. Furthermore, we would like to point
out that the reflective nature of the barrier, assumed in the
previous studies, is not essential for producing power law
distributions. Indeed, one only needs a positive lower barrier,
e.g., porous, absorbing or reflective one, which is a natural
condition since no physical object or socioeconomic one can approach
zero arbitrarily close without repelling from it or simply
disappearing. In many areas, objects of interest may not have a
strictly reflecting barrier, but rather a porous one, e.g., cities
may degenerate, bankruptcy protection may sometimes fail and a
company can be liquidated. In these cases, the power law effect
follows from the well-known queueing results on cycle maximum that
we briefly stated in Subsection \ref{ss:cycleMax}. This observation
presents a rigorous explanation for the previous study in
\cite{BS00} that argued heuristically how multiplicative processes
with absorbing barriers can result in power laws.

Here, we would like to point out that the stochastic difference
 equation ($M_{n+1}=J_nM_n+Q_n$) with random coefficients  is
 closely related to RMPs and is known to produce power law
 distributions. It appears that the
 first rigorous study of this process was done by Kesten in 1973
 \cite{Kesten73}; for later investigations of this model see \cite{HaanResnick89,Saporta04}
 and the references therein. In addition, we refer the reader to
 equations $(1.1)-(1.6)$ on p.126 of \cite{Goldie91} for other related
 stochastic recursions of multiplicative nature that produce power
 laws.

 Next, it is easy to see that RMBPs reduce to RMPs in the
special case when constant number of individuals are born in each
state of the modulating process. However,  our main result,
Theorem~\ref{theorem:rbp}, reveals a general asymptotic equivalence
between the power law exponent of a RMBP and the corresponding RMP.
In other words, Theorem \ref{theorem:rbp} discovers the asymptotic
insensitivity of the power law exponent on the conditional
distributions of the reflected branching process beyond their
conditional mean values. Furthermore, for the special case when the
modulating process is i.i.d., we sharpen the result on the
logarithmic asymptotics of Theorem \ref{theorem:rbp} to the exact
one in Theorem \ref{theorem:ExactMBP} by using the implicit renewal
theory of Goldie 1991 \cite{Goldie91}.

In some domains, e.g., the growth of living organisms, the objects
always grow (basically never shrink) up until a certain random time.
Huberman and Adamic \cite{AH99} also propose this model as an
explanation of the growth dynamics of the World Wide Web by arguing
that the observation time is an exponential random variable. This
notion has been revisited in \cite{R01} and generalized to a larger
family of random processes observed at an exponential random time
\cite{RH02}. In this regard, in Subsection~\ref{ss:smbp}, we study
randomly stopped modulated branching processes and show, under more
general conditions than the preceding studies, that the resulting
variables follow power laws.

In regard to the previously mentioned situations with absorbing
barriers,   we discuss MBP with an absorbing barrier in
Subsection~\ref{ss:porous}  and argue that it leads to power law
distributions as well.  We conjecture that these types of models can
be natural candidates for describing the bursts of requests at
popular Internet Web sites, often referred to as hotspots.

Based on our new model, we discuss two related phenomena: truncated
power laws and double Pareto distributions. We argue that one can
obtain a truncated power law distribution by adding an upper barrier
to RMBP, similarly as the truncated geometric distributions appear
in queueing theory, e.g., finite buffer $M/M/1$ queue. Furthermore,
by the duality of RMBP and queueing theory, we give two new natural
explanations of the origins of double Pareto distributions that have
been observed in practice. In the queueing context, it has been
shown that the tail of the queue length distribution exhibits
different decay rates in the heavy-traffic and large deviation
regime, respectively \cite{OBG06}; similar behavior of the queue
length distribution was attributed to the multiple time scale
arrivals in \cite{JL95}. We claim that the preceding two mechanisms,
when translated to the proportional growth context, provide natural
explanations of the double Pareto distributions.

Finally, we would like to mention that there might be other non
multiplicative mechanisms that result in power law distributions,
e.g., the randomly typing model used to explain the power law
distribution of frequencies of words in natural languages
\cite{Mi04} and the highly optimized tolerance studied in
\cite{CD99}.
Very recently, the new power law phenomenon in the situations where
jobs have to restart from the beginning after a failure was
discovered in \cite{FS05} and further studied in
\cite{SL06,Asmussen07}; equivalently in the communication context,
the  retransmission based protocols in data networks were shown to
almost invariably lead to power laws and, in general, heavy tails in
\cite{Jelen07ALOHA,PJ07RETRANS,Jelen07e2e,JT07t}. For a recent
survey on various mechanisms that result in power laws see
\cite{Mi04}.

The rest of the paper is organized as follows. After introducing the
modulated branching processes in Section \ref{s:model}, we study the
duality between the queueing theory and the multiplicative processes
with reflected barriers in Subsection \ref{ss:multiplicative} and
absorbing barriers in Subsection \ref{ss:cycleMax}, respectively.
Then, we present our main results in Section \ref{s:main} on the
logarithmic asymptotics of the stationary distribution of the
reflected modulated branching process and the corresponding
multiplicative one, which is followed by the study of the exact
asymptotics under the more restrictive conditions in Section
\ref{s:exact}. As further extensions, we discuss three related
models in Section \ref{s:relatedModels}, i.e., randomly stopped
processes in Subsection \ref{ss:randomStop}, modulated branching
processes with absorbing barriers in Subsection \ref{ss:porous} and
truncated power laws in Subsection \ref{ss:truncate}. In the end,
Section \ref{s:proofs} presents the majority of the technical proofs
that have been deferred from the preceding sections for increased
readability.

\section{Reflected Modulated Branching Processes}\label{s:model}
In this section we formally describe our model.
Let $\{J_n\}_{n > -\infty}$ be a stationary and ergodic modulating
process that takes values in positive integers.
    Define a family of independent, non-negative, integer-valued random variables
    $\{ B_n^i(j) \}, -\infty <i,j,n < \infty,$ which are independent of the modulating process $\{J_n\}$.
    In addition, for fixed $j$, variables $\{B(j), B_n^i(j) \}$ are
   identically distributed with $\mu(j)\eqdef
   \expect[B(j)]<\infty$.
   \begin{definition}\label{def:Zn}
    A Modulated Branching Process (MBP) $\{ Z_n \}_{n=0}^{\infty}$ is recursively defined by
        \begin{equation}
          Z_{n+1}\eqdef \sum_{i=1}^{Z_n} B_n^i(J_n) ,
          \label{eq:branching}
        \end{equation}
    where  the initial value $Z_0$ is a positive integer. For increased clarity,
    we may explicitly write $\{ Z_n^l \}$ when  $Z_0=l$.
   \end{definition}

    \begin{definition}\label{def:rmbp}
     For any $l \in \nat$ and an integer valued $\rbp_0$, a Reflected Modulated Branching Process (RMBP)
     $\{ \rbp_n \}_{n=0}^{\infty}$  is recursively defined as
         \begin{equation}
           \rbp_{n+1}\eqdef \max {\Bigg(} \; \sum_{i=1}^{\rbp_n} B_n^i(J_n) ,\; l \;{\Bigg)}. \;
            \label{eq:rbranching}
          \end{equation}
    \end{definition}
    \begin{remark}{\rm
These types of modulated branching processes with a reflecting
barrier appear to be new and, thus, the traditional methods for the
analysis of branching processes \cite{AV72} do not seem to directly
apply.}
\end{remark}
      \begin{remark}\label{remark:general}
       {\rm  A more general framework would be to define
         \begin{equation}\label{eq:general}
            Z_{n+1}=\int_{0}^{Z_n} B_n^t(J_n(t))d\nu(t),
            \end{equation}
         for any real measure $\nu$ and,  similarly,
         \begin{equation}
           \rbp_{n+1}=\max {\Bigg(} \int_{0}^{\rbp_n} B_n^t(J_n(t))d\nu(t) ,\; l \;{\Bigg)} ,
           \end{equation}
     where $l>0$ and $B_n^t(J_n(t))$ is $\nu$-measurable.
        We refrain from this generalization since it introduces additional technical difficulties without
        much new insight.
     }\end{remark}

Now, we present the basic limiting results on
        the convergence to stationarity of $Z_n$ and $\rbp_n$.

   \begin{lemma}
      \label{lemma:branch}
       If $\expect \log\mu(J_0)<0$, then a.s., we have
            \begin{equation*}
              \lim_{n\rightarrow \infty}Z_n=0.
             \end{equation*}
   \end{lemma}

 \begin{proof}
         For all $n\geq 1$,  let $W_n=Z_n / \Pi^0_{n-1}$, where $\Pi^0_n=\prod_{i=0}^{n}\mu(J_i)$.
         It is easy to check
         that $W_n$ is a positive
          martingale with respect to the filtration $\mathcal{F}_n=\sigma(J_i,Z_i,0\leq i \leq
          n-1)$. Hence, by the martingale convergence theorem (see Theorem 35.5. of \cite{PB05}),
          almost surely (a.s.) as $n\rightarrow
          \infty$,
            \begin{equation}\label{eq:converge}
             W_n\rightarrow W<\infty.
            \end{equation}
        Next, since $\{J_n\}$ is stationary and ergodic, so is $\{ \mu ( J_n) \}$, and
        therefore, a.s.,
           \[ \frac{\log \Pi^0_{n-1}}{n} =\frac{1}{n}\sum_{i=0}^{n-1} \log \mu(J_i) \rightarrow \expect \log\mu(J_0) <0
           \;\;\;  \text{as $n \rightarrow \infty$}.\]
         Thus, $\Pi^0_{n-1} \rightarrow 0 $  as $n \rightarrow \infty$,
         which, by recalling (\ref{eq:converge}) and $Z_n=W_n\Pi^0_{n-1}$,
        finishes the proof.
 \end{proof}

Next, let $Z_{-n}$ be the number of individuals at time $0$ in an
unrestricted  branching process that starts at time
    $-n$ with $l$ individuals; when needed for clarity, we will use the notation
    $Z_{-n}^{l}$ to explicitly indicate the initial state $l$.
 \begin{lemma} \label{lemma:zl}
  Assume $\expect \log\mu(J_0)<0$,
    then, for any a.s. finite initial condition $\rbp_0$,
    $\rbp_n$ converges in distribution to
                             \begin{equation*} \rbp \eqd \max_{n\geq 0} Z_{-n}.
                             \end{equation*}
  \end{lemma}

  \begin{proof}
      First, assume that $\rbp_0 =l$ and let $Z_n^k$ be the number of individuals at time $n$ in an unrestricted branching process that starts at time
        $k$ with $l$ individuals. Then, by stationarity of $\{J_n\}$, we have $Z_n^k\eqd Z_{k-n}$.  Clearly,
        \[    \rbp_1 = \max {\Bigg(}  \;  \sum_{i=1}^l B_1^i (J_1), l \; {\Bigg)}
                          \eqd   \max \{ \; Z_{-1}, Z_0 \;\},\]
                          and,
        by induction  and stationarity, it is easy to show
                    \begin{equation}\label{eq:repRBP}
                   \rbp_n \eqd \max( \; Z_{-n},Z_{-(n-1)},\cdots, Z_{-1},Z_0  \;) ,
                   \end{equation}
      where $Z_{-n},\cdots,Z_0 $ are defined on the same sequence $\{B_k^i(J_k)\}_{i\geq 1}$, ${-n \leq k \leq 0}$.
      Hence, by monotonicity, we obtain
               $$   \Pr[\rbp_n>x]\rightarrow \Pr[\rbp>x] \;\; \text{as}\; n \rightarrow \infty .$$

         Now, if $\rbp_{n}^{\rbp_{0}} $ is a process defined on the same sequence $\{B_n^i(J_n)\}$
         with the initial condition
         $\rbp_0  \geq l$, then, it is easy to see that
         $$ \rbp_{n}^{\rbp_{0}}   \geq \rbp_n \geq l, \;\;\; \text{for all $n$} ,$$
         implying
         \begin{equation}
         \label{eq:low}
          \Pr[\rbp_{n}^{\rbp_{0}}  > x] \geq \Pr[\rbp_n  > x].
          \end{equation}
         Next, if we define the stopping time $\tau$ to be the first time when $\rbp_{n}^{\rbp_{0}}  $ hits
         the boundary $l$, then the preceding monotonicity implies
         that $\rbp_n =\rbp_{n}^{\rbp_{0}} $ for all $n\geq \tau$.
          Using this observation, we obtain
         \begin{align}
         \label{eq:upper}
          \Pr[\rbp_{n}^{\rbp_{0}}   > x] &= \Pr[\rbp_{n}^{\rbp_{0}}   > x, \tau >n] +
                  \Pr[\rbp_{n}^{\rbp_{0}}   > x, \tau \leq n]  \nonumber \\
                            &=  \Pr[\rbp_{n}^{\rbp_{0}}   > x, \tau >n] + \Pr[\rbp_n > x, \tau \leq  n] \nonumber \\
                            &\leq  \Pr[\tau >n] + \Pr[\rbp_n > x].
         \end{align}
         Next,  by Lemma \ref{lemma:branch}, $\tau$ is a.s. finite and,
         thus, by (\ref{eq:low}) and  (\ref{eq:upper}),  we conclude
         $$  \lim_{n\rightarrow \infty} \Pr[\rbp_n  > x] = \lim_{n\rightarrow \infty} \Pr[\rbp_{n}^{\rbp_{0}}  > x]
                       = \Pr[\rbp > x]. $$
   \end{proof}

\subsection{Reflected Multiplicative Processes and Queueing Duality}\label{ss:multiplicative}
 Note that in the special case
$B_n^i(J_n)\equiv J_n$, reflected
  modulated branching processes reduce to reflected multiplicative processes with $J_n$ being integer valued.
  In general, using the definition in (\ref{eq:general}),  $J_n$ can be relaxed to
take any positive real values. Hence, in this subsection we assume
that  $\{\mdp_n\}_{n\geq 0}$ is a positive, real valued process.

\begin{definition}
  For $l>0$ and $M_0<\infty$, define a
   Reflected Multiplicative
    Process (RMP) as
      \begin{equation}
          \rmp_{n+1}=\max(\, \rmp_n \cdot \mdp_{n} \, , \, l \,), \;\;\;\; n\geq 0.
          \label{eq:rmp}
       \end{equation}
       \end{definition}

  The preceding RMP model was studied by Goldie in 1991 \cite{Goldie91};
  for later considerations of this model see \cite{SC97,LS96a1,LS96b,Ga99,GL05,Do01}.  Goldie \cite{Goldie91} also shows
   a direct connection (duality) between
  RMP and queuing theory in Section 5 of \cite{Goldie91} for the case when $\{J_n\}$ is an i.i.d. sequence.
  Here, we study this duality further in the
  generality of stationary and ergodic processes.

    Without loss of generality we can assume $l=1$, since we can always divide  (\ref{eq:rmp}) by $l$ and
    define $\rmp^1_n=\rmp_n/l$.
    Now, let $X_n=\log\mdp_n$ and $Q_n=\log\rmp_n $ with the standard conventions
    $\log0=-\infty$ and $e^{-\infty}=0$.
Then, for $l=1$, equation (\ref{eq:rmp}) is equivalent to
                 \begin{equation}
                  Q_{n+1}= \max(Q_n+X_{n},0),
                    \label{eq:queue}
            \end{equation}
   which is the workload (waiting-time) recursion in a single server (FIFO) queue.

    \begin{lemma} \label{lemma:Mrep}
                 If $\expect \log\mdp_n<0$, then $\rmp_n$ converges in distribution to an a.s. finite random
                 variable $\rmp$ that satisfies
                  \begin{equation} \rmp \eqd \sup_{n\geq 0}\Pi_n,
                  \end{equation}
                  where $\Pi_0=1, \Pi_n=\prod_{i=-n}^{-1}\mdp_i,n\geq 1$.
    \end{lemma}

   \begin{proof} By the classical result of Loynes \cite{Lo61}, $Q_n$, defined by (\ref{eq:queue}),
   converges in distribution to an a.s. finite stationary limit $Q$ if $\expect X_n=\expect \log\mdp_n
   <0$ and, furthermore,
               \[ Q\eqd \sup_{n\geq 0}S_n , \]
               where $S_0=0$ and $S_n=\sum_{i=-n}^{-1}X_i$.
  This implies the convergence in distribution of $\rmp_n$ to
  $$      \rmp \eqd e^{\sup_{n\geq 0}S_n}
                       = \sup_{n\geq 0}e^{S_n}
                       =\sup_{n\geq 0} \Pi_n.  $$
    \end{proof}

    The following theorem is a direct corollary of Theorem~1 in
    \cite{GW93}; see also Theorem 3.8 in \cite{Ch94} and, for a more
    recent presentation, we refer the reader to \cite{GOW04}.
    \begin{theorem}
    \label{theorem:multiplicative}
         Let $\{ J_n\}_{n\geq 1}$ be a stationary and ergodic sequence of positive random variables. If there exists a function $\Psi$ and positive
         constants $\alpha^{\ast}$ and $\varepsilon^{\ast}$ such that
         \begin{enumerate}
          \item[1)] $ n^{-1}\log \expect [(\Pi_n)^{\alpha}]\rightarrow \Psi(\alpha)$ as $n \rightarrow \infty$ for
          $\mid \alpha -\alpha^{\ast}\mid < \varepsilon^{\ast}$,
          \item[2)] $\Psi$ is finite and differentiable in a neighborhood of $\alpha^{\ast}$ with
          $\Psi(\alpha^{\ast})=0$, $\Psi'(\alpha^{\ast})>0$, and
          \item[3)] $\expect \left[(\Pi_n)^{\alpha^{\ast}+\varepsilon} \right]<\infty$, for $n\geq 1$ and some $\varepsilon >0$,
          \end{enumerate} then
            \begin{equation}    \lim_{x\rightarrow \infty} \frac{\log\Pr[\rmp> x]}{\log x}=-\alpha^{\ast}.
            \end{equation}
    \end{theorem}
\begin{remark}
{\rm
 We refer to conditions $1)-3)$ as the {\it polynomial G\"artner-Ellis
 conditions}. Note that condition $2)$ can be relaxed such that $\Psi$ is only differentiable at
 $\alpha^{\ast}$ and condition $3)$ can be weakened to $\varepsilon =0$ \cite{GW93}.
 Since conditions $2)$ and $3)$ are used for Theorem~\ref{theorem:rbp} in Section~\ref{s:main},  we keep the
 current form to provide a unified framework.
  Also, it is worth noting that the multiplicative process $\Pi_n$
without
 the reflective boundary would essentially follow the lognormal distribution, as it was
 recently observed in \cite{GL05} (this is similar to
 the fact that the unrestricted additive random walk is approximated well by Normal
 distribution). However, we would like to reemphasize that the
 lower boundary $l$ is not just a mathematical artifact, but a very
 natural condition since no physical object can approach zero arbitrarily
 close without either repelling (reflecting) from it or vanishing
 (absorbing); the absorbing boundary will be discussed in the following Subsection \ref{ss:cycleMax}.
 }\end{remark}
   Here, we illustrate the preceding theorem by the following examples. Assume
   that $\{A_n\},\{C_n\}$ are two mutually independent sequences, and let $J_n=e^{A_n-C_n}$.
    Then the quantity $Q_n\eqdef \log
   M_n$, where $M_n$ is defined in (\ref{eq:rmp}), satisfies
         \begin{equation}\label{eq:queue2} Q_{n+1}=(Q_n+A_n-C_n)^{+}. \end{equation}
   The first two examples assume
   that $\{A_n\},\{C_n\}$ are two i.i.d.
    sequences, the third example takes $\{J_n \}$ to be a Markov
    chain, and in the last example, $\{J_n \}$ is modulated by a
    Markov chain $\{X_n\}$.
         \begin{example}
         {\rm
      If $\{A_n\},\{C_n\}$ follow exponential distributions,
       $\Pr[C_n >x]=e^{-\mu x}$ , $\Pr[A_n >x]=e^{-\lambda x}$
      and $\lambda < \mu$,  then $Q_n$ represents the waiting time in a $M/M/1$ queue.  By
       Theorem 9.1 of \cite{As87}, the stationary waiting time in a $M/M/1$
       queue is distributed as
         \begin{equation*} \Pr[Q>x]=\frac{\lambda}{\mu}e^{-(\mu - \lambda)x} ,  \;\;\; x\geq 0,
         \end{equation*}
      which equivalently yields a power law distribution for $M$,
         \begin{equation*}  \Pr[M>x]=\Pr[Q>\log x]=\frac{\lambda}{\mu x^{\mu - \lambda}},\;\;\; x\geq 1
         \end{equation*}
         with power exponent $\alpha=\mu-\lambda$.
         }
\end{example}
   \begin{example}{\rm
      If $\{A_n\},\{C_n\}$ are two i.i.d Bernoulli processes with
      $\Pr[A_n=1]=1-\Pr[A_n=0]=p$,  $\Pr[C_n=1]=1-\Pr[C_n=0]=q$,
      $p<q$. Then, the elementary queueing/Markov chain theory shows
      that the stationary distribution of $Q_n$, as defined in
      (\ref{eq:queue2}), is geometric
      $\Pr[Q \geq j] = (1-\rho)\rho^{j}, j\geq 0$, where $\rho = p(1-q)/q(1-p) <1$. Therefore,
      \begin{align*}
        \Pr[M\geq x] &= \Pr[Q\geq \log x]= \rho^{\lfloor \log x \rfloor}, \;\;\; x\geq 1.
        \end{align*}
        Since $ \log x -1 < \lfloor \log x \rfloor \leq \log x $, it is easy to conclude that
          \begin{equation*}   \frac{1}{ x^{\log(1/\rho)} } \leq  \Pr[M\geq x] < \frac{1}{\rho x^{\log(1/\rho)} }.
          \end{equation*}

          }
 \end{example}

 \begin{example}{\rm
 If $\{J_n\}$ is a Markov chain taking values in a finite set $\Sigma$ and possessing an
 irreducible transition matrix $Q=\left( q(i,j)\right)_{ i,j \in \Sigma }$, then the function $\Psi$ defined in Theorem
 \ref{theorem:multiplicative} can be explicitly computed.  To this end, define matrix $Q_{\alpha}$ with elements
\begin{equation*}
  q_{\alpha}(i,j) = q(i,j) j^{\alpha}, \;\; i,j \in \Sigma.
\end{equation*}
 By Theorem 3.1.2 of
 \cite{DZ98}, we have as $n\to \infty$,
 \begin{equation*}
  n^{-1}\log \expect [(\Pi_n)^{\alpha}]\rightarrow \log \left({\rm
 dev}(Q_{\alpha}) \right),
 \end{equation*}
 where ${\rm dev}(Q_{\alpha})$ is the Perron-Frobenius eigenvalue of matrix $Q_{\alpha}$.
 To illustrate this result, we take $\Sigma=\{u, d\}$
 where $u =1/d>1$,
 and $q(d,u)=q, q(d,d)=1-q, q(u,d)=p, q(u,u)=1-p$ where $p>q$. It is easy to
 compute
\begin{equation*}
Q_{\alpha}= \left(
\begin{array}{cc}
 (1-p)u^{\alpha}  & pd^{\alpha} \\
      qu^{\alpha} & (1-q)d^{\alpha}  \\
\end{array}
 \right),
\end{equation*}
 and, by letting $\log \left({\rm  dev}(Q_{\alpha})\right)=0$,  we obtain
 \begin{equation*}
   \alpha^{\ast}=\frac{\log(1-q)-\log(1-p)}{\log u}.
 \end{equation*}
 }
\end{example}

\begin{example}[double Pareto] {\rm
\label{ex:dp} If $\{J_n \equiv J(X_n\}$ is modulated by a Markov
chain $X_n$, we argue that $\Pr[M>x]$ can have different asymptotic
decay rates over multiple time scales. This phenomenon was
investigated in \cite{JL95} in the queueing context and formulated
as Theorem 3 therein. To visualize this phenomena, we study the
following example. Consider a Markov process $X_n$ of two states
(say $\{1,2\}$) with transition probabilities $p_{12}=1/5000$,
$p_{21}=1/10$, and $\Pr[J(1)=1.2]=1-\Pr[J(1)=0.6]=0.5$,
$\Pr[J(2)=1.7]=1-\Pr[J(2)=0.25]=0.6$. The corresponding simulation
result for $5\times 10^7$ trials is presented in Figure
\ref{fig:doublePareto}.  We observe from this figure a double Pareto
distribution for $M$, which provides a new explanation to the
origins of double Pareto distributions as compared to the one in
\cite{Re01}.

 \begin{figure}[h]
\centering
\includegraphics[width=3.55in]{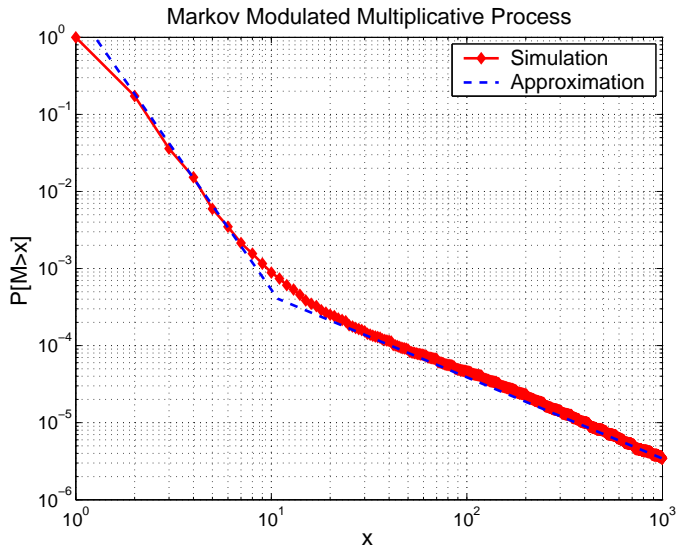}
  \caption{Illustration for Example~\ref{ex:dp} of the double Pareto distribution.}
\label{fig:doublePareto}
\end{figure}

}

\end{example}

\begin{remark} {\rm
For reasons of simplicity, we have chosen $\{J_n\}$ in all of the
preceding examples to be Markovian. However, Theorem
\ref{theorem:multiplicative} extends beyond the Markovian framework,
e.g., $\{J_n\}$ can be a semi-Markov process where the periods of
(sojourn) time that the process spends in a state are asymptotically
exponential but not necessarily memoryless. }
\end{remark}

\subsection{Multiplicative Processes with Absorbing Barriers and Cycle Maximum}\label{ss:cycleMax}
As briefly discussed in the introduction, we explained that the
reflective nature of the barrier is not essential for producing
power law distributions. Indeed, one only needs a positive lower
barrier, e.g., porous, absorbing or reflective one, which is a
natural condition since no physical objects or socioeconomic ones
can approach zero arbitrarily close without repelling from it or
simply disappearing. To illustrate the situations when the objects
can vanish, we name a few examples, e.g., cities may degenerate,
bankruptcy protection may sometimes fail and a company can be
liquidated. In these cases, the power law effect follows from the
well-known queueing result on cycle maximum that is
 stated in Theorem \ref{theorem:maximum} below. We also discuss in Subsection \ref{ss:porous} a more
 complicated situation when newly generated objects in the system can arrive/appear or leave/disappear.

  Following the notation from Chapter VIII of \cite{As87},  for a sequence of positive
  i.i.d. random variables $\{J, J_n \}_{n\geq 1}$,
   denote by $G_{+}$ the ladder height distribution of the random walk
          $\{ S_n=\sum_{i=1}^{n}\log J_i \}_{n \geq 1 }$ with
          $ \|G_{+}\| = \Pr[S_n\leq 0 \text{ for all } n \geq 1]$,
          and define the stopping time $\tau\eqdef \inf\{n: S_n\leq 0, n\geq 1 \}$
          with the corresponding cycle maximum $M_{\tau}\eqdef \sup \{\prod_{i=1}^{n}J_i: 1\leq n \leq \tau
          \}$;  here we assume, without loss of
          generality, that the absorbing barrier is equal to $1$.
     \begin{theorem}
      \label{theorem:maximum}
         If the sequence $\{ \log J_n \}_{n\geq 1}$ is  nonlattice,  satisfying
      $\expect[J^{\alpha^{\ast}}]=1$, $\alpha^{\ast}>0$ and
       $\left(\expect[J^{\alpha}]\right)'|_{\alpha=\alpha^{\ast}} <\infty$,
           then
                \begin{equation*}  \lim_{x \rightarrow \infty} \Pr[\rmp_{\tau}>x] x^{\alpha^{\ast}} =
               \frac{\left(1- \|G_{+}\| \right)\left(1-\expect\left[e^{-\alpha^{\ast}S_{\tau}} \right]\right) }{ \alpha^{\ast} \int_{0}^{\infty} xe^{\alpha^{\ast}x}G_{+}(dx)} .
               \end{equation*}
         \end{theorem}
         \begin{proof} The result is a direct consequence of Corollary 5.9 on p.~368
of \cite{As87}.
         \end{proof}

\section{Main Results}\label{s:main}
 This section presents our main results in Theorems \ref{theorem:rbp} and \ref{theorem:lessOne}.
 To avoid technical
 difficulties,  we assume $ \underline{\mu} \eqdef \inf_j \; \mu(j) > 0$. With a
  small abuse of notation,  as compared to the preceding Subsection \ref{ss:multiplicative}, we redefine here
   $\Pi_n=\prod_{i=-n}^{-1}\mu (\mdp_i), \, n\geq
  1, \; \Pi_0=l$ and $M=\sup_{n\geq 0}\Pi_n$.   In this paper we use the following
standard notation. For any two real functions $a(t)$ and $b(t)$, we
use $a(t)=o(b(t))$ to denote that $\lim_{t\rightarrow \infty}
a(t)/b(t)=0$, and $a(t)=O(b(t))$ to denote that
$\varlimsup_{t\rightarrow \infty} a(t)/b(t)<\infty$; when needed for
increased clarity, we may explicitly write $a(t)=o(b(t))$ as $t\to
\infty$.

 \begin{theorem}
 \label{theorem:rbp}
Assume that the process $\{\Pi_n\}$ satisfies the polynomial
G\"artner-Ellis conditions (conditions $1)-3)$ of Theorem
\ref{theorem:multiplicative}),  and $\sup_j \expect\left[e^{\theta
|B(j)-\mu(j)|}\right]< \infty$ for some $\theta >0$,
         then,
             \begin{equation}\label{eq:tworesults}
                 \lim_{x\rightarrow \infty} \frac{\log\Pr[\rbp>x]}{\log x}=
             \lim_{x\rightarrow \infty} \frac{\log\Pr[M>x]}{\log x}=-\alpha^{\ast}.
             \end{equation}
 \end{theorem}
\begin{remark}{\rm
Note that conditions 1) and 2) of
Theorem~\ref{theorem:multiplicative} imply that there exists $j$
such that $\mu(j)>1$, since otherwise we have
$\sup_{\alpha}\Psi(\alpha)\leq 0$, which would contradict
$\Psi(\alpha^{\ast})=0$ and $\Psi'(\alpha^{\ast})>0$ in condition
2). The following theorem covers the opposite situation when the
previous condition is not satisfied, i.e., $\sup_{j}\mu(j)<1$.
 }\end{remark}

  \begin{theorem}
  \label{theorem:lessOne}
        If $\sup_{j} \mu(j)<1$ and $\sup_j \expect\left[e^{\theta
|B(j)-\mu(j)|}\right]< \infty$ for some $\theta>0$, then,
 \begin{equation}
             \lim_{x \rightarrow \infty} \frac{\log\Pr[\rbp>x]}{\log x} = -\infty.
         \end{equation}

  \end{theorem}
\begin{remark}{\rm
  Informally speaking, these two theorems show that the alternating periods of contractions and expansions, e.g.,
   economic booms and recessions, are primarily responsible for the
appearance of power law distributions; in other words, if there are
no periods of expansions, i.e., the condition $\sup_{j} \mu(j)<1$ of
Theorem~\ref{theorem:lessOne} is satisfied, then $\Lambda$ has a
tail that is lighter than any power law distribution. Furthermore,
the first equality in (\ref{eq:tworesults}) of Theorem
\ref{theorem:rbp} reveals a general asymptotic equivalence between
the reflected modulated branching process and the corresponding
reflected multiplicative process, showing that the power law
exponent $\alpha^{\ast}$ is insensitive to the higher order
distributional properties of $B(j)$ beyond the conditional mean
$\mu(j)$. }
\end{remark}
\begin{remark}{\rm
A careful examination of the proofs reveals that the existence of a
uniform upper bound of the exponential moments for $|B(j)-\mu(j)|$
could possibly be relaxed to $\sup_j\expect\left[
|B(j)-\mu(j)|^{\alpha}\right]<\infty$ for $\alpha>\alpha^{\ast}$.
However, such an extension would considerably complicate the proofs.
Furthermore, in most practical applications the distributions of
$\{B(j)\}$ are typically very concentrated. For the preceding
reasons, we do not consider such extensions. }
\end{remark}
We present the \textbf{proofs} of Theorems~\ref{theorem:rbp} and
\ref{theorem:lessOne} in Subsection \ref{subsection:proofTh3}.

\section{Exact Asymptotics}\label{s:exact}
This section presents the exact asymptotic approximations of the
RMPs and RMBPs in the following two subsections, respectively.
\subsection{Exact Asymptotics of RMPs and the double Pareto phenomenon}
  The following two theorems essentially provide a new \rm{general} explanation of the measured
  double Pareto phenomenon (e.g., see \cite{Mi02,Re01}) since they
  rely on two universal statistical laws: the first one being based on the large deviation theory and the latter
  being implied by the central limit theorem.

The  theorems are direct translations from the corresponding
queueing theory results.  Theorem \ref{theorem:iidM} is based on the
large deviation result that studies the situation when $M$ is large,
and Theorem
  \ref{theorem:heavyTraffic} is derived from the heavy traffic approximation of a GI/GI/1 queue
  where we study the limiting behavior of
  a sequence of multiplicative processes with the multiplicative drift
  tending to one. These two theorems are basically corollaries of Theorem 5.2 in Chapter XIII
   and Theorem 7.1 in Chapter X of \cite{As87}, respectively.

   For a sequence of positive i.i.d. random variables $\{J, J_n \}_{n\geq 1}$,
   define $G_{+}$ to be the ladder height distribution of the random walk
          $\{ S_n=\sum_{i=1}^{n}\log J_i \}_{n \geq 1 }$ with
          $ \|G_{+}\| = \Pr[S_n\leq 0 \text{ for all } n \geq 1]$.
     \begin{theorem}
      \label{theorem:iidM}
         If the sequence $\{ \log J_n \}_{n\geq 1}$ is nonlattice, satisfying
           $\expect[J^{\alpha^{\ast}}]=1$,  $\alpha^{\ast}>0$, and
           $\left(\expect[J^{\alpha}]\right)'|_{\alpha = \alpha^{\ast}}<\infty$,
           then
                \begin{equation*}  \lim_{x \rightarrow \infty} \Pr[\rmp>x] x^{\alpha^{\ast}} =
               \frac{1- \|G_{+}\| }{ \alpha^{\ast} \int_{0}^{\infty} xe^{\alpha^{\ast}x}G_{+}(dx)} .
               \end{equation*}
         \end{theorem}
         \begin{proof} The result is a direct consequence of Theorem 5.3 in Chapter XIII
   of \cite{As87}.
         \end{proof}
         \begin{remark}{\rm
         If $S_n$ is lattice valued, see Remark 5.4 of Chapter XIII on p.~366 of \cite{As87}.
         }
         \end{remark}

Now, we study the limiting behavior of
  a sequence of multiplicative processes indexed by an integer $k$ where
  $J^{(k)}$, $S_n^{(k)}$ and $M^{(k)}$ are properly defined for all $k\geq 1$.
  \begin{theorem}
  \label{theorem:heavyTraffic}
     If $\left\{J^{(k)}, J_n^{(k)} \right\}_{n\geq 1}$ are positive and i.i.d. for each fixed $k$ with
      $m_k\eqdef \expect\left[\log J^{(k)}\right]$, $\sigma_k^2\eqdef \var\left[\log J^{(k)}\right]$,
      the random walks
          $\left\{ S_n^{(k)}=\sum_{i=1}^{n}\log J_i^{(k)} \right\}_{n \geq 1}$ satisfy $m_k<0$,
          $\lim_{k\rightarrow \infty }m_k=0$,
          $\underline{\lim}_{k\rightarrow \infty}\sigma_k^2>0$, and $\left(\log J^{(k)}\right)^2$ is
          uniformly integrable for all $k$,
          then, for $y\geq 1$,
          \begin{equation*}   \lim_{k\rightarrow \infty} \Pr\left[ \left(M^{(k)}\right)^{-m_k
          /  \sigma_k^2}>y \right]=1/y^2.
          \end{equation*}
  \end{theorem}
  \begin{proof}
        From Theorem 7.1 in Chapter X on p.~287 of \cite{As87}, we
        have, for $z\geq 0$,
      \[   \lim_{k\rightarrow \infty} \Pr\left[ - \frac{m_k}{\sigma_k^2  } \log M^{(k)}>z
           \right]= e^{-2z},
          \]
       which, by letting $z=\log y$, finishes the proof of Theorem \ref{theorem:heavyTraffic}.
  \end{proof}

\subsection{Exact Asymptotics of Reflected Branching Processes}
In this subsection, assuming that $\{J, J_n \}_{n\geq 1}$ are i.i.d.
and $\{\log \mu(J)\}$ is nonlattice, we will give an exact
asymptotics for RMBPs using the implicit renewal theorem of Goldie
(1991);  see Theorem 2.3 and Corollary 2.4 in \cite{Goldie91}. To
this end, let $\{B(j), B^i(j)\}_{i, j}$ be independent random
variables that are independent of $\{J, J_n\}$ and satisfy $B^i(j)
\eqd B(j)$.
\begin{theorem}\label{theorem:ExactMBP}
 If $\sup_{j} \expect \left[e^{\theta \mid B(j)-\mu(j) \mid
      } \right]<\infty$ for some $\theta>0$,
       $  \expect[\mu(J)^{\alpha^{\ast}}]=1
     $, $\alpha^{\ast}>0$ and
     $\expect[\mu(J)^{\alpha^{\ast}+\delta}]<\infty
     $ for some $\delta>0$, then,
      \begin{equation} \label{eq:ExactMBPresult}
       \lim_{ x\rightarrow \infty} \Pr[\Lambda ^{}>x]x^{\alpha^{\ast}}
         =\frac{ \expect
           \left[ {\left(\Lambda^{\ast}\right)}^{\alpha^{\ast}}-(\mu(J) \Lambda)^{\alpha^{\ast}}
           \right]}{\alpha^{\ast} \expect\left[\mu(J)^{\alpha^{\ast}} \log \mu(J)\right]},
         \end{equation}
     where $\Lambda^{\ast}\eqdef \max {\Bigg(} \; \sum_{i=1}^{\rbp^{}} B^i(J) ,\; l \;
     {\Bigg)}$ and $\Lambda$ is independent of $J$ and $\{B^i(j)\}_{i, j}$.
\end{theorem}

\begin{remark}{\rm
The preceding result is implicit because the constant on the right
hand side of equation  (\ref{eq:ExactMBPresult}) involves the value
of $\Lambda$, which is what we are trying to compute. In principle,
to derive the explicit exact asymptotics for RMBPs is a difficult
problem since the asymptotic constant depends on the behavior around
the boundary $l$. However, in the scaling region where the boundary
$l$ grows as well, albeit slowly, one can derive an explicit
asymptotic characterization.
}
\end{remark}

 In the following, similarly as in Theorem \ref{theorem:iidM},  we let $G_{+}$ be the ladder
height distribution of the nonlattice random walk
          $\{ S_n=\sum_{i=1}^{n}\log \mu(J_i) \}_{ n \geq 1}$ with
          $ \|G_{+}\| = \Pr[S_n\leq 0 \text{ for all } n \geq 1]<1
          $.
  \begin{theorem}
  \label{theorem:asympBrch}
      If $\sup_{k}\expect \left[e^{\theta  \mid B(k)-\mu(k) \mid
      } \right]<\infty$ for some $\theta>0$,
           $\underline{\mu}\eqdef \inf_j\mu(j)>0$, $  \expect[\mu(J)^{\alpha^{\ast}}]=1
     $ for some $\alpha^{\ast}>0$ with $\left(\expect[J^{\alpha}]\right)'|_{\alpha = \alpha^{\ast}}<\infty$,
         then, for any $\gamma>0$,
      \begin{equation*} {\large \lim_{\tiny \begin{array}{c}
  l_x\geq (\log x)^{3+\gamma} \\
  x\rightarrow \infty
\end{array}}} \Pr[\Lambda ^{l_x}/l_x>x]x^{\alpha^{\ast}}
         =\frac{1- \|G_{+}\| }{ \alpha^{\ast} \int_{0}^{\infty} xe^{\alpha^{\ast}x}G_{+}(dx)}.
         \end{equation*}

    \end{theorem}
\noindent The \textbf{proofs} of  Theorems \ref{theorem:ExactMBP}
and \ref{theorem:asympBrch} are presented in Subsection
\ref{subsection:theorem7}.  Here, we illustrate the exact asymptotics
of the reflected branching
    process with the following simulation example.

\begin{figure}[!th]
\label{fig:asymp}
\begin{center}
\includegraphics[height=3.1in]
{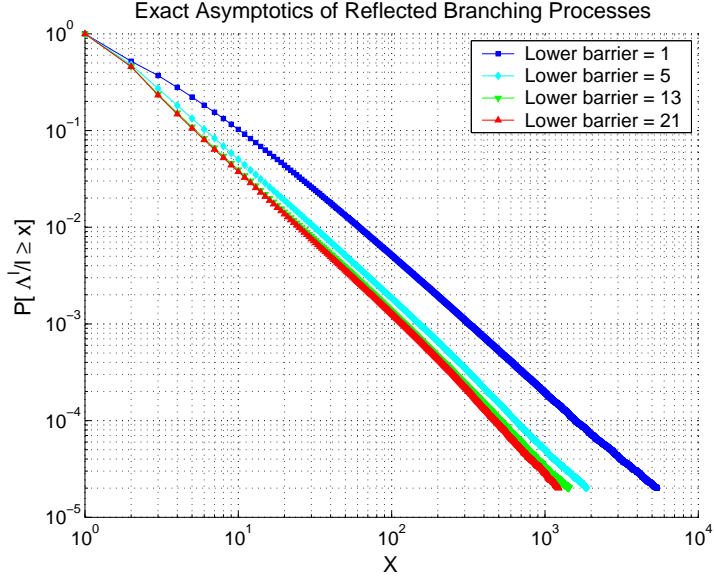} \caption{Simulation of $\Pr[\Lambda ^{l}/l\geq x]$
versus $x$ parameterized by $l$. }
\end{center}
\end{figure}

    \begin{example}{\rm
     Assume that $\{J_n\}_{n\geq 1}$ is a Bernoulli process with
      $\Pr[J_n=1]=0.4=1-\Pr[J_n=0]$ and the i.i.d. random variables $\{B_n^i(1)\}_{i\geq 1}$, $\{B_n^i(0)\}_{i\geq 1}$
       follow  Poisson
      distributions with means $1.5$ and $0.6$, respectively.
     The simulation results of $10^7$ samples, for $l=1,5,13,21$, are drawn in
     Figure \ref{fig:asymp}. From the figure we can clearly see that
     $\Pr[\Lambda ^{l}/l \geq x]$ approaches the limiting value very
     quickly, i.e., for $l=13$ and $l=21$, the plots of $\Pr[\Lambda^l/l\geq x]$ are basically
     indistinguishable.}
\end{example}

\section{Discussion of Related Models}\label{s:relatedModels}
Based on the study of reflected modulated branching processes, we
address three related models: randomly stopped processes, modulated
branching processes with absorbing barriers and truncated power
laws.

\subsection{Randomly Stopped Processes}\label{ss:randomStop}
In this subsection we discuss randomly stopped  multiplicative and
branching processes, respectively.
\subsubsection{Randomly Stopped Multiplicative Processes}
The following two theorems show that randomly stopped multiplicative
processes and reflected multiplicative processes are intimately
related and, to a certain extent, basically equivalent under more
restrictive conditions.  By following the approach of Chapter VIII
of \cite{As87}, we study the ladder heights of a multiplicative
process. For any RMP with i.i.d positive multiplicative increments,
the random variable $M$, as defined in Lemma \ref{lemma:Mrep},  can
be represented in terms of the ladder heights. To this end, define
$\Pi_n^0 \eqdef \prod_{i=0}^{n}J_i$ and the ladder height process
$\{H_i\}_{i\geq 1}$ of $\{ S_n=\sum_{i=1}^{n}\log J_i\}_{n \geq 1 }$
with  $ \|G_{+}\| = \Pr[S_n\leq 0 \text{ for all } n \geq 1]<1
          $ and $H_i^{e} \eqdef e^{H_i}$.

\begin{theorem}
 Suppose that $\{J, J_n\}_{n\geq 1}$ is a positive i.i.d. sequence with $\expect[\log J]<0$,
   then,
 \begin{equation} M \eqd \prod_{i=1}^{N}H_i^{e},
 \end{equation}
 where $N$ is independent of $\{H_i^{e}\}_{i\geq 1}$ and follows a geometric distribution
 $\Pr[N>n]=\|G_{+}\|^n$.
\end{theorem}
\begin{proof}
Based on the well-known Pollaczek-Khinchin representation (see
Chapter VIII of \cite{As87})
$$ \log M \eqd \sum_{i=1}^{N}H_i,
$$ where $N$ is independent of $\{ H_i \}$ with $\Pr[N>n]=\|G_{+}\|^n$,
it immediately follows that
$$\Pr[M>x]=\Pr \left[e^{\sum_{i=1}^{N}H_i}>x \right]=\Pr\left[\prod_{i=1}^{N}H_i^{e}>x \right].$$
\end{proof}

Conversely, we can prove that if the observation time has
exponential tail, the stopped process has a power law tail under
quite general conditions. Note that here we do not require $\{J_n\}$
to be an i.i.d. sequence.
\begin{theorem}\label{theorem:stopMul}
Let $N$ be an integer random variable independent of $\{J_n\}$ with
\begin{equation*}
 \lim_{x\rightarrow \infty}\frac{\log
\Pr[N>x]}{x} =-\lambda< 0.
\end{equation*}
For a positive ergodic and stationary process $\{J_n\}_{n\geq 0}$,
if $n^{-1}\log \expect\left[\left(\Pi_n^0\right)^{\alpha}\right] \to
\Psi(\alpha)<\infty$ as $n\to \infty$ in a neighborhood of
$\alpha^{\ast}>0$,  $\Psi(\alpha)$ is differentiable at
$\alpha^{\ast}$ with
          $\Psi(\alpha^{\ast})=\lambda$, $\Psi'(\alpha^{\ast})>0$ and
          $\expect \left[\left(\Pi_n^0\right)^{\alpha^{\ast}} \right]<\infty$ for $n\geq 1$,
 then,
 \begin{equation} \lim_{x\rightarrow \infty}\frac{\log \Pr\left[\Pi_{N}^0>x\right]}{\log x}
= -\alpha^{\ast}.
\end{equation}
\end{theorem}
\noindent The \textbf{proof} of Theorem \ref{theorem:stopMul} is
presented in Subsection \ref{ss:stopMulProof}.
\begin{remark}{\rm
This theorem generalizes the previous results from
\cite{AH99,R01,RH02} where only i.i.d. multiplicative increments are
considered.
 }\end{remark}

 Actually the following theorem shows that
randomly stopped multiplicative processes and reflected
multiplicative processes are basically equivalent under more
restrictive conditions. This equivalence is established using
classical results on $M/GI/1$ queue. In this regard,  we assume that
$\{ J_n\}_{n\geq 1}$ is an i.i.d. process, $\Pi_n^0$ is the
corresponding multiplicative process, $N$ is a geometric random
variable that is independent of $\Pi_n^0$ with $\Pr[N>n]=\rho^n,
0<\rho<1$, and $\bar{G}(t), t\ge 0$ is a complementary cumulative
distribution function.
\begin{theorem}\label{theorem:stopIID}
If a randomly stopped multiplicative process $\Pi_N^0$ satisfies
$$\Pr\left[\log J_1\leq
x\right]=\int_{0}^{x}\bar{G}(y)dy {\Big /}
\int_{0}^{\infty}\bar{G}(y)dy, \; x\ge 0$$ for some
$\bar{G}(\cdot)$,
   then, we can always construct a stationary RMP such that $\Pi_N^0 \eqd M$. Furthermore,
 if in addition
$\int_{0}^{\infty}e^{\alpha^{\ast}y}\bar{G}(y)dy=\rho^{-1}\int_{0}^{\infty}\bar{G}(y)dy$,
 $\int_{0}^{\infty}ye^{\alpha^{\ast}y}\bar{G}(y)dy<\infty$ for
$\alpha^{\ast}>0$, then,
         \begin{equation*}
       \lim_{x \rightarrow \infty} \Pr\left[M>x\right] x^{\alpha^{\ast}}=  \lim_{x \rightarrow \infty} \Pr\left[\Pi_{N}^0>x\right] x^{\alpha^{\ast}} =
    \frac{(1-\rho)\int_{0}^{\infty}\bar{G}(y)dy}{\alpha^{\ast}\rho \int_{0}^{\infty}ye^{\alpha^{\ast}y}\bar{G}(y)dy }.
       \end{equation*}
\end{theorem}

\noindent The \textbf{proof} of Theorem \ref{theorem:stopIID} is
presented in Subsection \ref{ss:stopMulProof}.

\subsubsection{Randomly Stopped Branching Processes}\label{ss:smbp}
In the following theorem, we extend Theorem~\ref{theorem:stopMul} of
the  preceding subsection to the context of randomly stopped
branching processes. Define $\Pi_n^0 \eqdef
\prod_{i=0}^{n}\mu(J_i)$.
\begin{theorem}\label{theorem:rsbp}
Suppose that $N$ is independent of $B_n^i(j)\geq 1$ for all $n,i,j$.
Then, under the same conditions as in Theorem~\ref{theorem:stopMul}
with $\expect \left[\left(\Pi_n^0\right)^{\alpha} \right]<\infty$
for $n\geq 1$ and $\Psi(\alpha)$ being differentiable  in a
neighborhood of $\alpha^{\ast}>0$, we obtain, for $\{Z_n\}_{n\geq
0}$ defined in (\ref{eq:branching}) with a bounded initial value
$Z_0<z_0<\infty$,
    \begin{equation*} \lim_{x\rightarrow \infty}\frac{\log \Pr[Z_{N}>x ]}{\log x}
=\lim_{x\rightarrow \infty}\frac{\log \Pr[\Pi^0_{N}>x ]}{\log x}=
-\alpha^{\ast}.
 \end{equation*}
\end{theorem}
\noindent The \textbf{proof} of this theorem is based on similar
arguments as in the proof of Theorem \ref{theorem:rbp}, and we defer
it to Subsection \ref{ss:stopMulProof}.

\subsection{Branching Processes with Absorbing Barriers}\label{ss:porous}
For many dynamic processes, e.g., city sizes, quite often when the
sizes of the  objects fall below a threshold, the whole object
disappears, e.g., urban decay. Therefore, it is natural to study
branching processes with absorbing barriers. As already discussed in
Subsection \ref{ss:cycleMax}, we know that a single object with an
absorbing barrier can result in power law distributions based on the
duality with the queueing cycle maximum. In this context, we can
also study a more complicated situation where the newly generated
objects can join the system and evolve together. This naturally
models the arrivals to popular Web sites (hotspots), since
information (news) is distributed according to a branching process,
e.g., user $A$ passes the information to $B$ and $C$; further $B$
may inform $D$, etc. Empirical examination shows that Web requests
follow power law distributions, e.g., see \cite{HA00,AH02}.

For a lower barrier $l>0$ and the modulated branching process
$\{Z^l_n\}_{n \geq 1}$ with $Z_0^l=l$ specified in
Definition~\ref{def:Zn}, define stopping time $P \eqdef \inf
\{n>0:Z^l_n\leq l \}$, where the modulating process $\{J, J_n\}$ is
 a sequence of i.i.d. random variables. This branching
 process, denoted by $Z_P$, vanishes completely after $P$; it is easy to prove that
$\expect[P]<\infty$ when $\expect[\log \mu(J_0)]<0$.

Let the arrivals $\{ A_n\}_{n>-\infty}$ be a sequence of i.i.d.
Poisson random variables with $\expect[A_n]=q>0$ that is independent
of other random variables. At time $n$, $A_n$ objects are generated
and join the system, each evolving according to an i.i.d. copy of
the modulated branching process $Z_P$.  Suppose that the system has
reached its stationarity with $N_n$ objects being in the system at
time $n$, and then, by Little's Law, $\expect[N_n] = q \expect[P]$.
Furthermore, assume that object $j$ observed at time $n=0$, if any,
is generated at time $(-P^r_j)$ with a size $Z^l_{-P^r_j}$,  where
the random variables $\{P^r_j\}$ are i.i.d. and follow the
equilibrium distribution of $P$. Then,
  the total size of all objects $Z_{s}$ observed at time $n=0$ in stationarity
  can be represented as
 \begin{equation*}
    Z_{s} = \sum_{j=1}^{N_0}Z^l_{-P^r_j}.
 \end{equation*}

Next, we show that $Z_s$ follows a power law. The proof of the
following theorem is essentially a corollary of Theorem
\ref{theorem:rbp}.

\begin{theorem}\label{theorem:rdbranch}
Under the conditions described in this subsection, if $\{\mu(J_n)\}$
satisfies $\inf_{j}\mu(j)>0$,  $\expect[\log
         \mu(J)]<0$,  $\expect[\mu(J)^{\alpha^{\ast}}]=1$ for some $\alpha^{\ast}>0$,
     $\expect[\mu(J)^{\alpha^{\ast}+\delta}]<\infty
     $ for some $\delta>0$,  and
         $\sup_j \expect\left[e^{\theta
|B(j)-\mu(j)|}\right]< \infty$ for some $\theta >0$, then,
 \begin{equation*}
\lim_{x \rightarrow \infty}\frac{\log \Pr[Z_s>x]}{\log x} =
-\alpha^{\ast}.
 \end{equation*}
\end{theorem}
\noindent The \textbf{proof} of this theorem is presented in
Subsection \ref{ss:absorbingProof}.

\subsection{Truncated Power Laws}\label{ss:truncate}
Truncated power laws have been observed empirically in many
practical situations where the studied objects have natural upper
boundaries. Here, we want to point out that by using the duality
between the modulated branching processes and the queueing theory,
one easily obtains truncated power laws when adding both a lower and
an upper barrier to the modulated branching process. To illustrate
this point, recall that $M/M/1/b$ queue with a finite buffer $b$
 results in a truncated geometric distribution for the number of
customers in the queue, and by the duality, it essentially follows
that in a proportional growth world with both a lower and an upper
barrier, truncated power laws can naturally arise, playing a similar
role as truncated exponential/geometric distributions do in an
additive world. Prior related work on this subject can be found in
\cite{SC97}.

\section{Proofs}\label{s:proofs}
\subsection{Proofs of Theorems \ref{theorem:rbp} and \ref{theorem:lessOne}}\label{subsection:proofTh3}
The proof of Theorem \ref{theorem:rbp}, composed of the upper bound
and the lower bound, and the proof of Theorem \ref{theorem:lessOne}
are presented in the following three subsections, respectively.
\subsubsection{Proof of Theorem \ref{theorem:rbp}: Upper Bound}  Since the proof is based on the change (increase)
of boundary $l$, we denote this dependence explicitly as $\Lambda^l
\equiv \Lambda$. According to Lemma \ref{lemma:zl}, the initial
value of $\{\Lambda_n\}$ has no impact on $\Lambda$ and, therefore,
in this subsection we simply assume that $\Lambda_0^l = l$. Before
stating the proof of the upper bound, we establish preliminary
Lemmas \ref{lemma:I_3}, \ref{lemma:B}, \ref{lemma:I_2} and
\ref{lemma:l}.

The first lemma shows that, most likely, the supremum  of $Z_n$
occurs for an index $n\leq x$.
  \begin{lemma}
  \label{lemma:I_3}
       For any $\beta>0$, the branching process $Z_n^l$ defined in (\ref{eq:branching})
       satisfies 
      $$ \sum_{n> x}^{\infty}\Pr\left[Z_n^{l} >x\right]= O\left( \frac{1}{x^{\beta}} \right ) \;\;
         \text{as $x\rightarrow \infty$}.$$
    \end{lemma}

   \begin{proof}
   Similarly as in the proof of Lemma~\ref{lemma:branch}, note
   that for
   $\Pi^0_{n-1}=\prod_{i=0}^{n-1}\mu(J_i)$, the stochastic process
   $W_n= Z_n^l / \Pi^0_{n-1}, n\geq 1$ is a martingale with respect to the
   filtration $\mathcal{F}_n=\sigma(J_i,Z_i,0\leq i \leq
          n-1)$ that satisfies $\expect[W_1]=1$.  Therefore, by recalling $\Pi_{n}=\prod_{i=-n}^{-1}\mu(J_i)$,
           we obtain, for any $\epsilon
          >0$,
   \begin{align}\label{eq:martgle}
      \Pr[Z_n^l>x] &= \Pr[W_n  \Pi_{n-1}^0 >x] = \Pr[(W_n e^{-\epsilon n}) (\Pi^0_{n-1} e^{\epsilon n}) > x] \nonumber \\
                   &\leq \Pr[W_n e^{-\epsilon n} > 1] + \Pr[\Pi_{n} e^{\epsilon n} >x] \nonumber \\
                   &\leq \expect[W_n e^{-\epsilon n}] + \Pr[\Pi_{n} e^{\epsilon n} >x] .
     \end{align}
     Next, by using the martingale property $\expect[W_n]=\expect[W_1]=1$,
     we derive
     \begin{equation}\label{eq:markov}
       \sum_{n> x}^{\infty} \expect[W_n e^{-\epsilon n}]=
      \sum_{n>x}^{\infty}e^{-\epsilon n}\leq \frac{e^{-\epsilon x}}{1-e^{-\epsilon }}
                 = O\left(\frac{1}{x^{ \beta }}\right) \;\; \text{as}\; x\rightarrow \infty.
     \end{equation}

     Then,  recalling conditions 1) and 2) of Theorem \ref{theorem:rbp} (or Theorem \ref{theorem:multiplicative}), we can choose $\delta , \epsilon >0$ small enough
      and $n_0$ large enough such that
      $ \Psi(\alpha^{\ast}-\delta)+2 \epsilon
      (\alpha^{\ast}-\delta)=-\zeta<0$ and
       $n^{-1}\log \expect\left[\Pi_n ^{(\alpha^{\ast}-\delta)}\right]<\Psi(\alpha^{\ast}-\delta)+
           \epsilon (\alpha^{\ast}-\delta)$ for $n>n_0$,  which implies, for $ x >n_0$,
      \begin{align} \label{eq:last}
        \sum_{n> x}^{\infty} \Pr[\Pi_{n} e^{\epsilon n} >x] &\leq \sum_{n> x}^{\infty}
                       \frac{\expect\left[\Pi_n ^{(\alpha^{\ast}-\delta)}\right]
                       e^{\epsilon (\alpha^{\ast}-\delta)
                       n}}{x^{(\alpha^{\ast}-\delta)}}
                       \leq   \sum_{n>x}^{\infty} \frac{e^{-\zeta n}}{x^{\alpha^{\ast}-\delta}}   \nonumber\\
                       &\leq \frac{e^{-\zeta x}}{(1-e^{-\zeta}) x^{\alpha^{\ast}- \delta} }
                       = O\left(\frac{1}{x^{\beta}} \right)\;\text{as}\; x \rightarrow
                       \infty.
      \end{align}
      Finally, by using  (\ref{eq:martgle}), (\ref{eq:markov}) and (\ref{eq:last}), we complete the proof.
      \end{proof}

The following lemma relates $\Lambda_n$ to the corresponding
multiplicative process.
    \begin{lemma}
    \label{lemma:B}
      Let $\epsilon>0$ and $\rbp_n^l$ be the reflected branching process, as defined in (\ref{eq:rbranching}),
      then,
      for $x\geq l$,
      $$ \Pr\left[\rbp_n^l>x\right]\leq \Pr\left[\max_{1\leq j \leq n} \mp_j (1+\epsilon)^j > x/l\right]+
         n\Pr\left[\mathcal{B}_0^{l,\epsilon}\right], $$
    where $\mp_j=\prod^{-1}_{i=-j} \mu (J_i)$ and $\mathcal{B}_n^{l,\epsilon}=
    \bigcup_{j\geq l}\{ \; \sum_{i=1}^j B_n^i(J_n)>j \mu (J_n) (1+\epsilon) \; \}$ for any integer $n$.
  \end{lemma}
 \begin{proof}
From (\ref{eq:repRBP}), we have
\begin{equation}\label{eq:repRBP2}
        \rbp_n^l \eqd \max\left( \; Z_{-n}^l,Z_{-(n-1)}^l,\cdots, Z_{-1}^l,Z_0^l  \;\right)
        .
\end{equation}
Next, let $Z_{-n}^{l}(k)$ be the branching process that starts at
time $-n$ with $l$ objects and is observed at time $k\geq -n$. Note
that $Z_{-i}^{l}(-i) = l$, $Z_{-i}^l(0)\equiv  Z_{-i}^l$ and
$$
         Z_{-j}^l = \sum_{i=1}^{Z_{-j}^{l}(-1)} B_{-1}^{i}(J_{-1})
$$
for $j\geq 1$. Now, by using the preceding observation,
(\ref{eq:repRBP2}) and $Z_0^l=l$, we derive, for $x\geq l$,
 \begin{align*}
   \Pr &\left[\rbp_n^l>x \right] \leq \Pr \left[ \max \left( \; Z_{-n}^l,Z_{-(n-1)}^l,\cdots, Z_{-1}^l, l  \;\right)>x,
   (\mathcal{B}_{-1}^{l,\epsilon})^C \right]+
                          \Pr \left[ \mathcal{B}_{-1}^{l,\epsilon} \right]\\
                          &\leq \Pr \left[ \max \left( \; Z_{-n}^l(-1)(1+\epsilon)\mu(J_{-1}),
                          \cdots, Z_{-2}^l(-1)(1+\epsilon)\mu(J_{-1}), l(1+\epsilon)\mu(J_{-1})  \; \right)>x,
                          (\mathcal{B}_{-1}^{l,\epsilon})^C \right] \nonumber\\
                          &\quad  +\Pr \left[ \mathcal{B}_{-1}^{l,\epsilon}\right]\nonumber\\
    &\leq \Pr \left[ \left\{ (1+\epsilon)\mu(J_{-1}) \max_{2\leq i \leq n} \left( Z_{i}^{l}(-1) \right)>x \right\} \; \bigcup \;
        \left\{ \mu(J_{-1})(1+\epsilon)>x/l \right\}  \right]
                      +\Pr \left[\mathcal{B}_{-1}^{l,\epsilon}
                      \right].
     \end{align*}
Then, intersecting with event $\mathcal{B}_{-2}^{l,\epsilon}$ and
using
$$
         Z_{-j}^l(-1) = \sum_{i=1}^{Z_{-j}^{l}(-2)} B_{-2}^{i}(J_{-2})
$$
for $j\geq 2$, one easily obtains
 \begin{align*}
   \Pr \left[\rbp_n^l>x \right]
    &\leq \Pr {\Bigg [} \left\{ (1+\epsilon)^2\mu(J_{-2}) \mu(J_{-1})
     \max_{3\leq i \leq n} \left( Z_{i}^{l}(-2) \right)>x \right\} \nonumber\\
     &\quad \quad \quad \bigcup \;
       {\Big \{} \max {\Big (} \; \mu (J_{-2})\mu (J_{-1}) (1+\epsilon)^2,
                    \;\mu (J_{-1}) (1+\epsilon) \; {\Big )}> x/l {\Big \}} {\Bigg
                    ]}\nonumber\\
     &\quad \quad +\Pr \left[\mathcal{B}_{-2}^{l,\epsilon} \right]+\Pr \left[\mathcal{B}_{-1}^{l,\epsilon}
                      \right],
     \end{align*}
which, by continuing the induction and using $\Pr
\left[\mathcal{B}_{i}^{l,\epsilon}\right]=\Pr\left[\mathcal{B}_{0}^{l,\epsilon}
\right]$ for all $i$,  finishes the proof.
 \end{proof}

Now, we show that the ``error'' event $\mathcal{B}_0^{l,\epsilon}$
in the preceding lemma has a negligible probability for large $l$
relative to any power law distribution.
 \begin{lemma}
  \label{lemma:I_2}
   If $\sup_j \expect\left[e^{\theta |B(j)-\mu(j)|}\right]<\infty$, $\theta>0$ and $\underline{\mu}=\inf_j \mu(j)>0$,
    then, by setting $l_x=\lfloor x^{\delta} \rfloor$, $0<\delta<1$ in the definition of
      $\mathcal{B}_0^{l_x,\epsilon}$ in Lemma~\ref{lemma:B},
    we obtain,
      for any $\beta>0$,
      \begin{equation*}
      \Pr\left[\mathcal{B}_0^{l_x,\epsilon}\right]=O\left( \frac{1}{x^{\beta}}
      \right) \;\;\; \text{as} \;\;\; x\rightarrow \infty.
      \end{equation*}
    \end{lemma}
 \begin{proof}
   First,   we derive
      \begin{align}
        {P}(n)  &\eqdef \Pr \left[ \sum_{i=1}^{n} B_0^i(J_0) > \mu(J_0)(1+\epsilon)n
        \right]\nonumber\\
        &         \leq \Pr \left[ \sum_{i=1}^{n} \left(B_0^i(J_0)-\mu(J_0)\right) >
\epsilon \underline{\mu} n \right] \nonumber\\
&\leq \left(\expect\left[ e^{\zeta\left(B(J) - \mu(J)\right)
}\right]\right)^n e^{-\zeta \epsilon \underline{\mu} n},  \quad
\zeta>0, \nonumber
     \end{align}
which, by using the elementary inequality
     $e^x\leq 1+x+x^2e^{|x|}/2$, $x\in \real$ and setting $x=B(J) - \mu(J)$,  yields
\begin{align}
 {P}(n)  & \leq \left( 1+ \frac{\zeta^2}{2} \expect\left[\left( B(J) - \mu(J)
\right)^2 e^{\zeta |B(J) - \mu(J)|}\right] \right)^{n} e^{-\zeta
\epsilon \underline{\mu} n}. \nonumber
\end{align}

For any $\eta>0$ and large enough $n$ such that $\zeta =\eta \log
n/( \epsilon \underline{\mu} n )< \theta$, the assumption
$\sup_{j}\expect\left[e^{\theta |B(j) - \mu(j)|}\right] < \infty$
implies
$$C\eqdef \expect\left[\left( B(J) - \mu(J) \right)^2 e^{\zeta |B(J) -
\mu(J)|}\right] < \infty,$$ which yields
\begin{align}\label{eq:P2}
 {P}(n)  & \leq \left( 1+ \frac{C(\eta \log n)^2}{2(n \epsilon \underline{\mu})^2}  \right)^{n} n^{-\eta}=O\left(\frac{1}{n^{\eta}} \right).
\end{align}

    Therefore,  choosing $\eta= 1+\beta/\delta$ in (\ref{eq:P2}),  we obtain, for $l_x=\lfloor x^{\delta} \rfloor,
    0<\delta<1$ and $\beta>0$,
    as $x \rightarrow \infty$,
     \begin{equation*}
      \Pr\left[\mathcal{B}_0^{l_x,\epsilon}\right] \leq   \sum_{i=l_x}^{\infty}{P}(n)
          \leq
          O\left(  \sum_{n= \lfloor x^{\delta} \rfloor}^{\infty} \frac{1}{n^{\eta}}  \right)
          = O\left(\frac{1}{x^{\beta}}\right).
     \end{equation*}
\end{proof}

The following lemma allows us to increase the lower barrier in order
to prove the upper bound.
  \begin{lemma}
  \label{lemma:l}
     Assume that $\rbp_n^{l_1}$ and $\rbp_n^{l_2}$ are defined on the same sequence $\{B_n^j(J_n)\}$ with initial
     conditions $l_1$ and $l_2$, respectively. If $l_1\geq l_2$, then, for
      all $n\geq 0$,
         $$  \rbp_n^{l_1}\geq \rbp_n^{l_2}. $$
  \end{lemma}
  \begin{proof}
  The result holds trivially for $n=0$. Now we prove the result using induction.
  Suppose that it is true for all $0\leq k \leq n$, and for $k=n+1$,
    \begin{align*}
     \rbp_{n+1}^{l_1} &= \max {\Bigg(} \; \sum_{i=1}^{\rbp_n^{l_1}} B_n^i(J_n) ,\; l_1 \;
     {\Bigg)}
     \geq \max {\Bigg(} \; \sum_{i=1}^{\rbp_n^{l_2}} B_n^i(J_n) ,\; l_2 \;{\Bigg)}= \rbp_{n+1}^{l_2},
    \end{align*}
  which implies that the lemma is true for all $n\geq 0$.
  \end{proof}

 Now, we are ready to complete the proof of the upper bound.

 \begin{proof}[Proof of the upper bound of Theorem \ref{theorem:rbp}:]
    Choosing $l_x=\lfloor x^{\epsilon} \rfloor \geq l, 0
    <\epsilon<1$,
    using Lemma~\ref{lemma:l} and then Lemma \ref{lemma:B}, we derive
     \begin{align}\label{eq:upMBP}
      \Pr\left[\rbp^{l} > x \right] &= \Pr \left[ \sup_{j\geq 1}Z_{-j}^{l}>x \right]
      \leq \Pr \left[ \rbp_{\lfloor x \rfloor}^{l}>x  \right] + \Pr \left [\sup_{j > x}Z_{-j}^{l}>x \right]
               \nonumber   \\
                          &\leq \Pr\left[\rbp_{\lfloor x \rfloor}^{l_x}>x \right] +\sum_{j>x} \Pr\left[Z_j^{l}>x \right] \nonumber\\
                          &\leq \Pr \left[\sup_{j \geq 1} \mp_j (1+\epsilon)^j > x^{1-\epsilon} \right]+
                          x\Pr \left[\mathcal{B}_0^{l_x,\epsilon} \right]   + \sum_{j>x} \Pr\left[Z_j^{l}>x \right]  \nonumber\\
                          &\eqdef I_1(x) +I_2(x) +I_3(x).
     \end{align}

  Now, define a new process $\{ \mu^{\epsilon}(J_n)= \mu(J_n)
(1+\epsilon)\}_{n\geq 1}$ and
     $\Pi_n^{\epsilon} = \prod_{i=-n}^{-1} \mu^{\epsilon}(J_i)  $. Then, for
     $\epsilon$ small enough, we have
      \begin{enumerate}
      \item[1)] $ n^{-1}\log \expect (\Pi_n^{\epsilon})^{\alpha}\rightarrow \Psi^{\epsilon}(\alpha)=
      \Psi(\alpha)+\alpha \log(1+\epsilon)$ as
          $n \rightarrow \infty$ for
          $\mid \alpha -\alpha^{\ast}\mid < \varepsilon^{\ast}$,
      \item[2)]  $\Psi^{\epsilon}$ is finite in a neighborhood of $\alpha^{\ast}_{\epsilon}$,
       $\alpha_{\epsilon}^{\ast} < \alpha^{\ast}$, and differentiable at
          $\alpha^{\ast}_{\epsilon}$ with
          $\Psi(\alpha_{\epsilon}^{\ast})+\alpha_{\epsilon}^{\ast}
          \log(1+\epsilon)=0$,
          $\Psi'(\alpha_{\epsilon}^{\ast})>0$, and
      \item[3)] $\expect \left[(\Pi_n^{\epsilon})^{\alpha^{\ast}_{\epsilon}} \right]<\infty$ for $n\geq
      1$.
      \end{enumerate}
     Therefore, by Theorem \ref{theorem:multiplicative}, we obtain
    \begin{equation} \label{eq:I3-1}
     \lim_{x\rightarrow \infty}
     \frac{\log \Pr[ \sup_{i \geq 1 } \mp_i  (1+\epsilon)^i> x^{1-\epsilon}]}{\log x} =
     -(1-\epsilon)\alpha^{\ast}_{\epsilon},
     \end{equation}
 which, in conjunction with Lemma \ref{lemma:I_3} and Lemma
    \ref{lemma:I_2}, yields
 \begin{equation} \label{eq:I3-2}
     I_2(x)+I_3(x)=o(I_1(x)).
     \end{equation}
  Then, combining (\ref{eq:upMBP}), (\ref{eq:I3-1}) and (\ref{eq:I3-2})
  yields
    \begin{align*}
      \frac{\log \Pr[\rbp^{l} > x ]}{\log x} &\leq \frac{\log \left( (1+o(1))I_1(x) \right)}{\log x}
       \longrightarrow  -(1-\epsilon)\alpha^{\ast}_{\epsilon} \;\;\;\text{as}\;  x\rightarrow \infty.
    \end{align*}
    Since $\Psi^{\epsilon}(\alpha)$ is continuous in a neighborhood
    of $\alpha^{\ast}$ in both $\alpha$ and $\epsilon$, we derive
    \begin{equation*} \lim_{\epsilon \rightarrow 0}\alpha^{\ast}_{\epsilon} =
    \alpha^{\ast}, \end{equation*}
    implying,
    \begin{equation}\label{eq:upinequality}
      \varlimsup_{x \rightarrow \infty}\frac{\log \Pr[\rbp > x ]}{\log x}
    \leq
    -\alpha^{\ast}.\end{equation}
 \end{proof}

\subsubsection{Proof of Theorem \ref{theorem:rbp}: Lower Bound}
In order to prove the lower bound, we need to establish
the following three lemmas.  Specifically,  Corollary
\ref{cor:copies} allows us to obtain a lower bound for $\Lambda$
while, maybe somewhat counterintuitively, increasing the lower
barrier $l$.
\begin{lemma}\label{lemma:2lambda}
  Let $\{ \rbp_n^{y_1} \}$ and $\{ \rbp_n^{y_2} \}$ be defined on the same modulating sequence
 $\{ J_n\}_{n\geq 0}$ and independent random  variables $\left\{B_n^{i,1}(j), B_n^{i,2}(j)\right\}$
 with $B_n^{i,k}(j)$ identically distributed for fixed $j$. Then,
  \begin{equation*}
  \rbp_n^{y_1+y_2}\leqd  \rbp_n^{y_1}+\rbp_n^{y_2} ,
  \end{equation*}
  where $\rbp_n^{y_1}$ and $\rbp_n^{y_2}$ are conditionally
  independent given $\{J_n\}_{n \geq 0}$.
 \end{lemma}
 \begin{proof}
 We use induction to prove this lemma. Starting with $n=1$,  we obtain
       \begin{align*}
        \rbp_1^{y_1+y_2} &= \max {\Bigg (} \; \sum_{i=1}^{y_1+y_2} B_0^i(J_0), \; y_1+y_2 \; {\Bigg )}
        \eqd \max {\Bigg (} \; \sum_{i=1}^{y_1} B_0^{i,1}(J_0)+
                         \sum_{i=1}^{y_2} B_0^{i,2}(J_0), \; y_1+y_2 \; {\Bigg )}  \\
                         &\leq  \max {\Bigg (} \; \sum_{i=1}^{y_1} B_0^{i,1}(J_0), \; y_1\; {\Bigg )}
                         + \max {\Bigg (} \; \sum_{i=1}^{y_2} B_0^{i,2}(J_0), \; y_2 \; {\Bigg
                         )}\\
                         &\eqd \rbp_1^{y_1}+\rbp_1^{y_2},
       \end{align*}
       since  for any $x_1$,$x_2$,$y_1$,$y_2$,
  $$ \max (\; x_1+x_2\; , \; y_1+y_2 \;) \leq \max ( \; x_1\; ,y_1 \; ) + \max ( \; x_2\; ,y_2 \; ) . $$
   The proof is completed by induction in $n$,
       \begin{align*}
        \rbp_{n+1}^{y_1+y_2} &= \max {\Bigg (} \; \sum_{i=1}^{\rbp_n^{y_1+y_2}} B_1^i(J_n), \; y_1+y_2 \; {\Bigg )}
        \leqd \max {\Bigg (}\; \sum_{i=1}^{\rbp_n^{y_1}+\rbp_n^{y_2}}
                            B_1^i(J_n), \; y_1+y_2 \; {\Bigg )} \\
                         &\leqd \rbp_{n+1}^{y_1}+\rbp_{n+1}^{y_2}.
       \end{align*}
 \end{proof}
 Next, a straightforward application of the preceding lemma yields
 the following corollary.
  \begin{corollary}\label{cor:copies}
    If $\{ \rbp_{n,j}^{1} \}_{1\leq j \leq y}$  are  conditionally i.i.d copies of $\rbp_n^1$
     given $\{J_i\}_{1\leq i \leq n}$, then,
      $$  \rbp_n^{y} \leqd  \sum_{j=1}^{y}\rbp_{n,j}^{1} \;. $$
  \end{corollary}

Now, we basically establish that the supremum of $\Pi_i$ occurs most
likely for small indexes $i\leq h\log x$.
  \begin{lemma}
  \label{lemma:lowerI_2}
   Assume that condition 1) of Theorem \ref{theorem:rbp} is satisfied,
   then,  for $0\leq \epsilon<1$ and any $\beta>0$, there exists $h>0$ such that,  when $x \rightarrow \infty$,
       \[
         \Pr \left[ \sup_{i > h \log x} \mp_i  (1-\epsilon)^i> x \right]=
         O\left( \frac{1}{ x^{\beta} } \right).
     \]
  \end{lemma}
  \begin{proof}
  Using condition 1) of Theorem \ref{theorem:rbp}, we can
  choose   $0<\alpha< \alpha^{\ast}$ with $n^{-1}\log \expect[\Pi_n^{\alpha}] \rightarrow \Psi(\alpha)<0$
   and $n_0$  large enough, such that $\expect \left[ \mp_n^{\alpha} \right] <
   \zeta^n$,
  $0<\zeta<1$, $n>n_0$.  Thus, for  $h=-\beta/ \log
   \zeta>0$ and $x>e^{n_0/h}$,
  \begin{align*}
        \Pr \left[ \sup_{i > h \log x} \mp_i  (1-\epsilon)^i> x \right]
          &\leq \sum_{i> h \log x}^{\infty} \Pr \left[ \mp_i > x \right]
          \leq \sum_{i> h \log x}^{\infty} \frac{ \expect \left[ \mp_i^{\alpha} \right]
          }{ x^{\alpha} } \leq  \sum_{i> h \log x}^{\infty} \frac{ \zeta^i }{  x^{\alpha}}
          =O\left( \frac{1}{ x^{ \beta} }\right)   .
  \end{align*}
  \end{proof}

Finally, the last lemma shows that $\sum_{i=1}^j B_n^i(J_n)$ can not
deviate too much from $j\mu(J_n)$ for large $j$.
   \begin{lemma}
  \label{lemma:lowerI_3}
  For $0<\delta, \epsilon<1$ and $\mathcal{C}_n^{l,\epsilon} \eqdef
  \bigcup_{j\geq l}\{ \; \sum_{i=1}^j B_n^i(J_n) < j \mu (J_n) (1-\epsilon) \; \}$,
  we obtain, for any $\beta>0$,
  \begin{equation*}
         \Pr\left[\mathcal{C}_0^{\lfloor x^{\delta} \rfloor,\epsilon}\right] = O\left( \frac{1}{ x^{\beta} } \right)  .
 \end{equation*}
  \end{lemma}
\begin{proof}
The proof of this lemma is basically the same as Lemma
\ref{lemma:I_2}. Observe
     \begin{align*}
        {P}(n)  &\eqdef \Pr \left[ \sum_{i=1}^{n} B_0^i(J_0) < \mu(J_0)(1-\epsilon)n
        \right]\nonumber\\
        & \leq \Pr \left[ \sum_{i=1}^{n} \left(B_0^i(J_0)-\mu(J_0) \right) < -\epsilon {\mu(J_0)}n \right] \\
         & \leq \Pr \left[ \sum_{i=1}^{n} \left(\mu(J_0) - B_0^i(J_0) \right) > \epsilon \underline{\mu} n
         \right].\nonumber
     \end{align*}
 Then, by using a similar argument as in deriving (\ref{eq:P2}), we
can prove, for any $\beta>0$,
 \begin{equation*}
      \Pr\left[\mathcal{C}_0^{\lfloor x^{\delta} \rfloor,\epsilon}\right]= O\left(\frac{1}{x^{\beta}}\right).
     \end{equation*}
 \end{proof}

Next, we can complete the proof of the lower bound of Theorem
\ref{theorem:rbp}.

  \begin{proof}[Proof of the lower bound of Theorem \ref{theorem:rbp}:]
    First, using Corollary \ref{cor:copies}, we obtain, for any integer $y\geq 1$,
     \begin{align}\label{align:low}
       \Pr[\rbp_n^l > x] & \geq \Pr[\rbp_n^1 >x]  = \frac{y  \Pr[\rbp_n^1 >x]}{y}
       \geq \frac{\Pr[\sum_{j=1}^{y} \rbp_{n,j}^1 > y x]}{y} \geq  \frac{\Pr[\rbp_n^y
       >yx]}{y}.
    \end{align}

    Now, using (\ref{eq:repRBP}), similarly as in the proof of Lemma
       \ref{lemma:B}, for
$0<\epsilon <1$ and
$\mathcal{C}_n^{l,\epsilon}$ defined in Lemma \ref{lemma:lowerI_3},
we
    derive
     \begin{align}\label{eq:I123}
       \Pr[\rbp_n^y >yx]&\geq \Pr \left[ \max_{0\leq i \leq n} \left( Z_{-i}^{yx} \right)>yx ,
       \left(  \bigcup_{i=-n}^{-1}\mathcal{C}_{i}^{y,\, \epsilon}  \right)^C\right] \nonumber\\
&\geq \Pr \left[ \sup_{1 \leq i \leq n}\mp_i (1-\epsilon)^{i}> x,
\left(  \bigcup_{i=-n}^{-1}\mathcal{C}_{i}^{y,\, \epsilon}
\right)^C
\right] \nonumber\\
        &\geq \Pr \left[ \sup_{1 \leq i \leq n} \mp_i  (1-\epsilon)^i> x \right] - n \Pr[\mathcal{C}_0^{y,\, \epsilon}] \nonumber\\
        &\geq \Pr \left[ \sup_{i \geq 1} \mp_i  (1-\epsilon)^i> x
        \right]   -  \Pr \left[ \sup_{i > n} \mp_i  (1-\epsilon)^i> x
               \right] - n \Pr\left[\mathcal{C}_0^{y,\, \epsilon}\right] \nonumber\\
              &\eqdef I_1(x) - I_2(x)-I_3(x);
    \end{align}
    note that $\{I_j(x)\}_{1\leq j \leq 3}$ here are different from those in (\ref{eq:upMBP}).

     Next, similarly as in the proof of the upper bound, define a new process
      $\{ \mu_{\epsilon}(J_n)= \mu(J_n) (1-\epsilon)\}_{n\geq 1}$ and
     let $\Pi_n^{\epsilon} = \prod_{i=-n}^{-1} \mu_{\epsilon}(J_i)  $.
     Then, for  $\epsilon$ small enough, we have
      \begin{enumerate}
      \item[1)] $ n^{-1}\log \expect (\Pi_n^{\epsilon})^{\alpha}\rightarrow \Psi^{\epsilon}(\alpha)=
      \Psi(\alpha)+\alpha \log(1-\epsilon)$ as
          $n \rightarrow \infty$ for
          $\mid \alpha -\alpha^{\ast}\mid < \varepsilon^{\ast}$,
      \item[2)]  $\Psi^{\epsilon}(\alpha)$ is finite in a neighborhood of $\alpha^{\ast}_{\epsilon}$,
       $\alpha_{\epsilon}^{\ast} > \alpha^{\ast}$ and differentiable at
          $\alpha^{\ast}_{\epsilon}$ with
          $\Psi(\alpha_{\epsilon}^{\ast})+\alpha_{\epsilon}^{\ast}
          \log(1-\epsilon)=0$,  $\Psi'(\alpha_{\epsilon}^{\ast})>0$, and
      \item[3)] $\expect \left[(\Pi_n^{\epsilon})^{\alpha^{\ast}_{\epsilon}} \right]<\infty$ for $n\geq
      1$.
      \end{enumerate}
     Therefore, by Theorem \ref{theorem:multiplicative}, we
     obtain
    \begin{equation}\label{eq:I1}
     \lim_{x\rightarrow \infty}
     \frac{\log \Pr[ \sup_{i \geq 1 } \mp_i  (1-\epsilon)^i> x]}{\log x} =
     -\alpha^{\ast}_{\epsilon}.
     \end{equation}
     Now, by setting $y=\lfloor x^{\delta}
    \rfloor, 0<\delta<1$,
$n=\lfloor x \rfloor$ in (\ref{align:low}), (\ref{eq:I123}),
 and using Lemmas \ref{lemma:lowerI_2} and \ref{lemma:lowerI_3},  it
 is easy to see that
  \begin{equation}
   I_2(x)+I_3(x)=o(I_1(x)), \nonumber
  \end{equation}
  which, by (\ref{align:low}) and
(\ref{eq:I123}), yields
\begin{align}\label{eq:I123A}
        \log \Pr[\rbp > x] &\geq \log \Pr[\rbp_{n}^l > x] \nonumber\\
        &\geq  \log (I_1(x)-I_2(x)-I_3(x)) -  \delta \log
                          x\nonumber\\
        & = \log ((1-o(1))I_1(x)) - \delta \log x. \nonumber
     \end{align}
From the preceding inequality and (\ref{eq:I1}), we obtain
     \begin{equation}\label{eq:lowerinequality}
      \varliminf_{x\rightarrow \infty} \frac{\log \Pr[\rbp > x]}{\log x}
     \geq -\alpha^{\ast}_{\epsilon}-\delta.
      \end{equation}
   Since $\Psi^{\epsilon}(\alpha)$ is continuous in a neighborhood of $\alpha^{\ast}$ in both
   $\alpha$ and $\epsilon$, we have
   $
     \lim_{\epsilon \rightarrow 0}\alpha^{\ast}_{\epsilon} = \alpha^{\ast}
   $.
Then, passing $\epsilon, \delta \to 0$ in (\ref{eq:lowerinequality})
 completes the proof
of the lower bound, which,  in conjunction with
(\ref{eq:upinequality}), finishes the proof of
Theorem~\ref{theorem:rbp}.
   \end{proof}

 \subsubsection{Proof of Theorem \ref{theorem:lessOne}}
    \begin{proof}
Using the same arguments as in deriving (\ref{eq:upMBP}) in the
proof of the upper bound  of Theorem \ref{theorem:rbp}, we obtain,
    for $l_x=\lfloor x \rfloor \geq l$ and
    $0<\epsilon<1$,
      \begin{align}\label{eq:lessoneproof}
      \Pr\left[\rbp > x \right]
       & \leq \Pr \left[ \rbp_{\lfloor x \rfloor}^{l}>x \right] + \Pr \left [\sup_{j >  x }Z_{-j}^{l}>x \right]
               \nonumber   \\
                          &\leq \Pr\left[\rbp_{\lfloor x \rfloor}^{l_x}>x \right] +\sum_{j>x } \Pr\left[Z_j^{l}>x \right] \nonumber\\
                          &\leq \Pr \left[\sup_{j \geq 1} \mp_j (1+\epsilon)^j > 1 \right]+
                          x\Pr\left[\mathcal{B}_0^{l_x, \,\epsilon} \right]
                                 + \sum_{j>x }^{\infty} \Pr\left[Z_j^{l}>x\right] \nonumber \\
                          &\eqdef I_1(x) +I_2(x) +I_3(x).
     \end{align}
    Recalling $\Pi_j=\prod_{i=-1}^{-j}\mu(J_i)$ and noting $\sup_j \mu (j) < 1$, we can choose $\epsilon>0$
      such that $\sup_j \mu (j)(1+\epsilon)<1$, which implies $I_1(x)=0$.
     And,
 using Lemma \ref{lemma:I_2},  we obtain $I_2(x)=O\left(x^{-\beta}\right)$
for all $\beta>0$.

 Next, using similar arguments as in deriving (\ref{eq:martgle}) in the proof of Lemma
  \ref{lemma:I_3}, we obtain, for $\epsilon>0$ and $j\geq 1$,
   \begin{align*}
      \Pr\left[Z_j^l>x\right] &\leq \expect\left[W_j e^{-\epsilon j}\right] + \Pr\left[\Pi_{j} e^{\epsilon j}
      >x\right],
     \end{align*}
     which, by recalling $\sup_{j}\mu(j)<1$ and choosing $\epsilon$ small enough such that $\Pr[\Pi_{j} e^{\epsilon j}
      >x]=0$ for $x>1$, yields,
     \begin{equation*}
      I_3(x)\leq \sum_{j> \lfloor x \rfloor}^{\infty} \expect[W_j e^{-\epsilon j}]=
      \sum_{j>\lfloor x \rfloor}^{\infty}e^{-\epsilon j}=O\left(e^{-\epsilon
      x}\right).
     \end{equation*}

 Finally, combining (\ref{eq:lessoneproof}) and the
bounds on $I_1(x)$, $I_2(x)$ and $I_3(x)$ finishes the proof.
  \end{proof}

\subsection{Proofs of Theorems \ref{theorem:ExactMBP} and \ref{theorem:asympBrch}}\label{subsection:theorem7}
\subsubsection{Proof of Theorem \ref{theorem:ExactMBP}}
In order to prove Theorem \ref{theorem:ExactMBP},  we first derive
the following lemma.
 \begin{lemma}\label{lemma:ExactMBP}
Under the assumptions of Theorem \ref{theorem:ExactMBP}, there
exists  $\gamma >0 $,  such that
\begin{equation}
\expect \left[ \left|  \sum_{i=1}^{\Lambda} \left(B^i(J) -
\mu(J)\right) \right|^{\alpha^{\ast}+\gamma} \right]< \infty.
\nonumber
\end{equation}
 \end{lemma}
 \begin{proof}
 We observe that, for $x>0$,
\begin{align}
\Pr \left[ \left|  \sum_{i=1}^{\Lambda} \left(B^i(J) - \mu(J)\right)
\right| > x\right] &= \Pr \left[ \sum_{i=1}^{\Lambda} \left(B^i(J) -
\mu(J)\right)  > x\right]+\Pr \left[
 \mu(J) -\sum_{i=1}^{\Lambda} \left(B^i(J)  \right)  >
x\right] \nonumber\\
&\eqdef I_1 + I_2. \nonumber
\end{align}

 We choose $\gamma$ and $\delta$ such that $0<\gamma < \alpha^{\ast}$ and $\gamma< \delta< \alpha^{\ast}$.  To evaluate $I_1$,
 for $0< \epsilon< (\alpha^{\ast}-\delta)/2$,  we set
$0<\beta\eqdef (\delta+\epsilon)/(\alpha^{\ast} - \epsilon)<1$, and
obtain
\begin{align}\label{equation:I1bExactMBP}
I_1 &\leq \Pr\left[\Lambda > x^{1+\beta}\right] + \Pr\left[
\sum_{i=1}^{\Lambda} \left(B^i(J) - \mu(J)\right)  > x,  \Lambda
\leq x^{1+\beta}\right] ,
\end{align}
which, by recalling that Theorem \ref{theorem:rbp} implies
$\Pr[\Lambda>x]=O\left(x^{-(\alpha^{\ast} - \epsilon)} \right)$,
results in
\begin{equation}\label{equation:I1aExactMBP}
\Pr\left[\Lambda
> x^{1+\beta}\right]=O\left( x^{-(\alpha^{\ast}+\delta)}
\right).
\end{equation}
 Now, we  study the second probability on the right-hand side
 of (\ref{equation:I1bExactMBP}). Using the fact that $J$ is independent of $\Lambda$ and applying Chernoff bound,  we
 obtain
 for $\zeta>0$,
\begin{align}\label{equation:I1bExactMBP}
\Pr \left[ \sum_{i=1}^{\Lambda} \left(B^i(J) - \mu(J)\right)  > x,
\Lambda \leq x^{1+\beta}\right] &= \sum_{n=l}^{\lfloor x^{1+\beta}
\rfloor}\Pr[\Lambda = n] \Pr\left[ \sum_{i=1}^{n}\left(B^i(J) -
\mu(J)\right)    > x \right]\nonumber\\
&\leq  \sum_{n=l}^{\lfloor x^{1+\beta} \rfloor}\Pr[\Lambda = n]
e^{-\zeta x}\left(\expect\left[e^{\zeta\left(B(J) - \mu(J)\right)}
\right]\right)^n .
\end{align}
Then,  setting $t=\zeta\left(B(J) - \mu(J)\right)$ in
(\ref{equation:I1bExactMBP}), using $e^t\leq 1+t+t^2e^{|t|}/2$,
$t\in \real$ and observing that $\expect[ B(J) - \mu(J)]=0$,
(\ref{equation:I1bExactMBP}) is further upper bounded by
\begin{align}
&\sum_{n=l}^{\lfloor x^{1+\beta} \rfloor}\Pr[\Lambda = n] \left( 1+
\frac{\zeta^2}{2} \expect\left[\left( B(J) - \mu(J) \right)^2
e^{\zeta |B(J) - \mu(J)|}\right] \right)^n e^{-\zeta x}.\nonumber
\end{align}
Next,  $\sup_{j}\expect\left[e^{\theta |B(j) - \mu(j)|}\right]
<\infty$ implies $C\eqdef \expect\left[\left( B(J) - \mu(J)
\right)^2 e^{\zeta |B(J) - \mu(J)|}\right]<\infty$, $\zeta<\theta$.
Hence,  for $x$ large,  we have $\zeta = (\alpha^{\ast}+\delta)\log
x/x< \theta$, implying that (\ref{equation:I1bExactMBP}) is bounded
by
\begin{align}\label{equation:I11ExactMBP}
&\sum_{n=1}^{\lfloor x^{1+\beta} \rfloor}\Pr[\Lambda = n] \left( 1+
\frac{C((\alpha^{\ast}+\delta)\log x)^2}{2x^2} \right)^{x^{1+\beta}}
x^{- (\alpha^{\ast}+\delta)}=O\left( x^{-(\alpha^{\ast}+\delta)}
\right)
\end{align}
since $\beta<1$.  Combining (\ref{equation:I1aExactMBP}),
(\ref{equation:I1bExactMBP})
 and (\ref{equation:I11ExactMBP}) proves
\begin{equation}\label{equation:I1ExactMBP}
I_1=O\left( x^{-(\alpha^{\ast}+\delta)} \right).
\end{equation}

By using the same approach as in proving
(\ref{equation:I1ExactMBP}), we can also show
\begin{equation}
I_2=O\left( x^{-(\alpha^{\ast}+\delta)} \right), \nonumber
\end{equation}
implying
\begin{align}
\Pr \left[ \left|  \sum_{i=1}^{\Lambda} \left(B^i(J) - \mu(J)\right)
\right| > x\right] =O\left( x^{-(\alpha^{\ast}+\delta)} \right).
\nonumber
\end{align}

Therefore, since $\delta>\gamma$,
 \begin{align}
 \expect \left[ \left|
\sum_{i=1}^{\Lambda} \left(B^i(J) - \mu(J)\right)
\right|^{\alpha^{\ast}+\gamma} \right]&=\int_{0}^{\infty} \Pr \left[
\left|  \sum_{i=1}^{\Lambda} \left(B^i(J) - \mu(J)\right)
\right|^{\alpha^{\ast}+\gamma}
> x\right]dx \nonumber\\
&=O\left(\int_{0}^{\infty}
x^{-\frac{\alpha^{\ast}+\delta}{\alpha^{\ast}+\gamma}}dx\right)<\infty,
\nonumber
\end{align}
which finishes the proof of the lemma.
 \end{proof}

Now, we proceed with proving Theorem \ref{theorem:ExactMBP}.

\begin{proof}[Proof of Theorem \ref{theorem:ExactMBP}]
The proof is based on Corollary 2.4 in \cite{Goldie91}, for which it
is sufficient to show
\begin{align}\label{eq:checkExamMBP}
I \eqdef  \expect \left[ \left| \left( \max {\Bigg(} \;
\sum_{i=1}^{\rbp^{}} B^i(J) ,\; l \;
     {\Bigg)} \right)^{\alpha^{\ast}}-(\mu(J) \Lambda)^{\alpha^{\ast}}
     \right| \right] < \infty.
\end{align}

In order to prove the preceding inequality,  we will use the
following elementary inequality (see equation (9.27) in
\cite{Goldie91}), for $x, y \geq 0$,
\begin{equation}\label{eq:twoCase}
\left| x^{\alpha} - y^{\alpha} \right| \leq \left\{ \begin{aligned}
         &|x-y|^{\alpha}, & &0<{\alpha}\leq 1\\
         &\alpha |x-y|\left( x^{{\alpha}-1} + y^{{\alpha}-1}\right), &
         &1<{\alpha}<\infty.
                          \end{aligned} \right.
\end{equation}

First, we prove the case when $0<\alpha^{\ast}\leq 1$. Using the
fact that
\begin{align}\label{eq:I1cExactProof}
\left(\max \left( \; \sum_{i=1}^{\rbp^{}} B^i(J) ,\; l \;
     \right)\right)^{\alpha^{\ast}} &= l^{\alpha^{\ast}} \ind\left( \sum_{i=1}^{\rbp^{}} B^i(J)\leq l
     \right) + \left(\sum_{i=1}^{\rbp^{}} B^i(J)\right)^{\alpha^{\ast}}\nonumber\\
     &\quad\quad
      - \left(\sum_{i=1}^{\rbp^{}} B^i(J)\right)^{\alpha^{\ast}}\ind\left( \sum_{i=1}^{\rbp^{}} B^i(J)\leq l
      \right),
\end{align}
we obtain, by  using (\ref{eq:twoCase}) and recalling Lemma
\ref{lemma:ExactMBP},
\begin{align}\label{}
I &\leq  2l^{\alpha^{\ast}}+\expect \left[ \left|
\left(\sum_{i=1}^{\rbp^{}}
     B^i(J)\right)^{\alpha^{\ast}}-(\mu(J) \Lambda)^{\alpha^{\ast}}
     \right|\right]\nonumber\\
     &\leq 2l^{\alpha^{\ast}}+\expect \left[ \left| \sum_{i=1}^{\rbp^{}}
     \left(B^i(J)-\mu(J)\right)
     \right|^{\alpha^{\ast}}\right] <\infty.
\end{align}

Next, we prove the case when $\alpha^{\ast}>1$. Applying
(\ref{eq:twoCase}) and (\ref{eq:I1cExactProof}) in
(\ref{eq:checkExamMBP}), we obtain
\begin{align}
I &\leq  2l^{\alpha^{\ast}}+\expect \left[ \left|
\left(\sum_{i=1}^{\rbp^{}}
     B^i(J)\right)^{\alpha^{\ast}}-(\mu(J) \Lambda)^{\alpha^{\ast}}
     \right|\right]\nonumber\\
     &\leq 2l^{\alpha^{\ast}}+ \alpha^{\ast}\expect \left[ \left| \sum_{i=1}^{\rbp^{}}
     \left(B^i(J)-\mu(J)\right)
     \right| \left| \sum_{i=1}^{\rbp^{}}
     B^i(J) \right|^{\alpha^{\ast}-1}  \right]
     \nonumber\\
     &\quad \quad \quad +\alpha^{\ast} \expect \left[ \left| \sum_{i=1}^{\rbp^{}}
     \left(B^i(J)-\mu(J)\right)
     \right| \left|\Lambda \mu(J) \right|^{\alpha^{\ast}-1}
     \right]\nonumber\\
     &\eqdef 2l^{\alpha^{\ast}}+I_1+I_2.
\end{align}
 For $I_1$, we use H\"older's inequality to obtain, for
 $0<\epsilon<1/2$,
 \begin{align}\label{eq:I1aExactProof}
I_1 &\leq \alpha^{\ast}\expect \left[ \left| \sum_{i=1}^{\rbp^{}}
     \left(B^i(J)-\mu(J)\right)
     \right|^{\frac{\alpha^{\ast}}{1-\epsilon}} \right]^{\frac{1-\epsilon}{\alpha^{\ast}}}
     \expect \left[ \left( \sum_{i=1}^{\rbp^{}}
     B^i(J) \right)^{\frac{\left(\alpha^{\ast}-1\right)\alpha^{\ast}}
     {\alpha^{\ast}+\epsilon-1}}\right]^{\frac{\alpha^{\ast}+\epsilon-1}{\alpha^{\ast}}}\nonumber\\
&\leq \alpha^{\ast}\expect \left[ \left| \sum_{i=1}^{\rbp^{}}
     \left(B^i(J)-\mu(J)\right)
     \right|^{\frac{\alpha^{\ast}}{1-\epsilon}} \right]^{\frac{1-\epsilon}{\alpha^{\ast}}}
     \expect \left[ \Lambda ^{\frac{\left(\alpha^{\ast}-1\right)\alpha^{\ast}}
     {\alpha^{\ast}+\epsilon-1}}\right]^{\frac{\alpha^{\ast}+\epsilon-1}{\alpha^{\ast}}}
 \end{align}
 where the last inequality uses the fact that
 $$
\sum_{i=1}^{\rbp^{}} B^i(J) \leqd \Lambda.
 $$
Now, Theorem \ref{theorem:rbp} implies that
$$\expect \left[ \Lambda^{\frac{\left(\alpha^{\ast}-1\right)\alpha^{\ast}}
     {\alpha^{\ast}+\epsilon-1}}\right]<\infty,$$
which, in combination with  Lemma \ref{lemma:ExactMBP} and
(\ref{eq:I1aExactProof}),  results in
\begin{equation}\label{I1ExactProof}
I_1<\infty.
\end{equation}

Using the same argument as in proving (\ref{eq:I1aExactProof}), we
obtain, by noting that $\Lambda$ and $\mu(J)$ are independent from
each other,
 \begin{align}
I_2 &\leq \alpha^{\ast}\expect \left[ \left| \sum_{i=1}^{\rbp^{}}
     \left(B^i(J)-\mu(J)\right)
     \right|^{\frac{\alpha^{\ast}}{1-\epsilon}} \right]^{\frac{1-\epsilon}{\alpha^{\ast}}}
     \expect \left[ \Lambda^{\frac{\left(\alpha^{\ast}-1\right)\alpha^{\ast}}
     {\alpha^{\ast}+\epsilon-1}}\right]^{\frac{\alpha^{\ast}+\epsilon-1}{\alpha^{\ast}}}
     \expect \left[ \mu(J)^{\frac{\left(\alpha^{\ast}-1\right)\alpha^{\ast}}
     {\alpha^{\ast}+\epsilon-1}}\right]^{\frac{\alpha^{\ast}+\epsilon-1}{\alpha^{\ast}}}< \infty, \nonumber
 \end{align}
 which, in conjunction with (\ref{I1ExactProof}), proves
 (\ref{eq:checkExamMBP}) and finishes the proof of Theorem \ref{theorem:ExactMBP}.
\end{proof}

\subsubsection{Proof of Theorem \ref{theorem:asympBrch}}\label{ss:asympBrchProof}
 The proof of this
theorem relies on the following lemmas; the first one is based on
Theorem~3.7.1 of \cite{DZ98}.
    \begin{lemma}\label{lemma:ModDev}
        $\{X(j), X_i(j)\}_{i,j \in \mathbb{Z}}$ are zero mean independent random
        variables that are identically distributed for fixed $j$.
        Fix a sequence ${a_n} \rightarrow
        0$ such that $na_n \rightarrow \infty$ as $n \rightarrow
        \infty$.  If $\sup_j \expect \left[e^{\theta X(j)} \right]<\infty$
        for $\theta>0$,  then, there exist $n_0, h>0$, such that for all
        $n>n_0$ and any random variable  $J \in \mathbb{Z}$,
        \[ \Pr \left[ \, \sum_{i=1}^{n}X_i(J)> \sqrt{\frac{n}{a_n}} \, \right] \leq e^{-\frac{h}{a_n}} .
        \]
      \end{lemma}

\begin{proof}
  For $0<\beta< \theta$, define $\varphi_{J}(\omega) = \expect \left[e^{\omega X(J)} \mid J
   \right]$ and use Taylor  expansion to derive
   \[  \varphi_{J}(\omega)= \varphi_{J}(0)  + \varphi_{J}'(0) \omega +
   \frac{\varphi_{J}''(\zeta)}{2!} \omega^{2},  \;  0< \zeta <
   \omega \leq \beta.
   \]
     Noting that $\varphi_{J}(0) = 1$,  $\varphi_{J}'(0) =
   0$ and  $$\varphi_{J}''(\zeta)=\expect \left[ X(J)^2 e^{\zeta X(J)} \mid J
   \right] \leq \sup_{j}\expect \left[ X(j)^2 e^{\beta {X(j)}}
   \right]\eqdef K_{\beta}<\infty,$$ we obtain, for $0<\omega\leq \beta$,
   \[
          \varphi_{J}(\omega)\leq 1+ K_{\beta} \omega^2,
   \]
which implies,
    \begin{align*}
       \Pr \left[ \, \sum_{i=1}^{n}X_i(J)> \sqrt{\frac{n}{a_n}} \, \right]
       &\leq e^{-\omega \sqrt{\frac{n}{a_n}}} \expect\left[ \left(\varphi_{J}(\omega) \right)^n
       \right]
       \leq e^{-\omega \sqrt{\frac{n}{a_n}}} (1+ K_{\beta}
       \omega^2 )^n
       \leq e^{-\omega \sqrt{\frac{n}{a_n}}}e^{n K_{\beta}\omega^2}.
   \end{align*}
  Since there exists $n_0$
   such that $\beta> 1/(2K_{\beta}\sqrt{na_n})$ for all $n>n_0$, we
   can choose $\omega = 1/(2K_{\beta}\sqrt{na_n})$, which implies,
   for $n>n_0$,
  \begin{align*}
      \Pr \left[ \, \sum_{i=1}^{n}X_i(J)> \sqrt{\frac{n}{a_n}} \, \right]
         &\leq e^{-\frac{1}{4K_{\beta}a_n}
         }=e^{-\frac{h}{a_n}},
  \end{align*}
where $h= 1/(4K_{\beta})>0$.
\end{proof}

\begin{lemma}\label{lemma:rare}
For any $l \in \nat$, define
$$\mathscr{D}_n^{l} \eqdef
    \bigcup_{k \geq l} \left\{ \; \sum_{i=1}^k B_n^i(J_n)>
        k \mu (J_n) \left(1+  \underline{\mu}^{-1}l^{-\frac{1}{3}} \right) \;
     \right\}$$
      and
      $$\mathscr{E}_n^{l}\eqdef \bigcup_{k \geq l}\left\{ \;
      \sum_{i=1}^k
B_n^i(J_n) < k \mu (J_n) \left(1- \underline{\mu}^{-1}
l^{-\frac{1}{3}}\; \right) \; \right\}.$$
 If $l \geq (\log x)^{3+ \gamma }$, $\gamma >0$, then, under the conditions of Theorem
 \ref{theorem:asympBrch},  we obtain, for any $\beta>0$,
     $$ \Pr\left[\mathscr{D}_n^{l}\right]=O\left( \frac{1}{ x^{\beta }} \right)
      \;\; \text{and} \;\; \Pr\left[\mathscr{E}_n^{l}\right]=O\left( \frac{1}{ x^{\beta }} \right). $$
\end{lemma}
\begin{proof}
Defining $a_n=n^{-1/3}$ and observing that $n a_n$ is monotonically
increasing in $n$, we obtain
  \begin{align*}
     \Pr\left[\mathscr{D}_n^{l}\right] &\leq \sum_{k=l}^{\infty}
      \Pr \left[\; \sum_{i=1}^k B_n^i(J_n) > k \mu (J_n) \left(1+\frac{1}{\underline{\mu}} \sqrt{\frac{1}{l a_l}}\;\right) \right]\\
            &\leq \sum_{k=l}^{\infty}
            \Pr \left[\; \sum_{i=1}^k B_n^i(J_n) > k \mu (J_n)
                \left(1+\frac{1}{\underline{\mu}} \sqrt{\frac{1}{k a_k}}\; \right) \right]\\
&\leq \sum_{k=l}^{\infty} \Pr \left[\; \sum_{i=1}^k
\left(B_n^i(J_n)-\mu(J_n)\right)  > \sqrt{\frac{k}{a_k}} \;\right],
\end{align*}
which,   by applying Lemma \ref{lemma:ModDev}, yields
\begin{align*}
     \Pr\left[\mathscr{D}_n^{l}\right]
  &\leq \sum_{k=l}^{\infty} e^{-\frac{h}{ a_k}}
   \leq \sum_{k\geq (\log x)^{3+\gamma}}^{\infty} e^{-hk^{1/3}}=O\left( \frac{1}{
x^{\beta}} \right) \; \text{as $x \to \infty$.}
\end{align*}

By the same argument,
     \begin{align*}
     \Pr\left[\mathscr{E}_n^{l}\right]&\leq \sum_{k=l}^{\infty}
      \Pr \left[\; \sum_{i=1}^k B_n^i(J_n) < k \mu (J_n) \left(1-\frac{1}{\underline{\mu}} \sqrt{\frac{1}{l\cdot a_l}}\;\right) \right]\\
            &\leq \sum_{k=l}^{\infty} \Pr \left[\; \sum_{i=1}^k
\left( \mu(J_n)-B_n^i(J_n) \right)  > \sqrt{\frac{k}{a_k}}
\;\right]\\
&=O\left( \frac{1}{ x^{\beta}}  \right) \; \text{as $x \to \infty$}.
     \end{align*}
\end{proof}

Following the proof of Lemma \ref{lemma:I_3} with minor
modifications, we can prove the following stronger result.
\begin{lemma}
  \label{lemma:asymp_3}
   For any $\beta>0$,  there exists $h>0$ such that the branching process defined in (\ref{eq:branching}) satisfies
      $$ \sum_{n\geq h\log x}^{\infty}\Pr\left[Z_n^{l} >x\right]= O\left( \frac{1}{x^{\beta}} \right ) \;\;
         \text{as $x\rightarrow \infty$}.$$
    \end{lemma}

Now, we proceed with the proof of Theorem \ref{theorem:asympBrch}.

\begin{proof}[Proof (of Theorem \ref{theorem:asympBrch}):]
 First, we establish the {upper bound}.  Setting $\epsilon = \underline{\mu}^{-1}
l^{-1/3}$ in Lemma \ref{lemma:B}, we
   obtain
     $$ \Pr \left[\rbp_n^l> lx \right]\leq \Pr\left[\max_{1\leq j \leq n} \mp_j
     \left(1+\underline{\mu}^{-1}
l^{-\frac{1}{3}} \right)^j > x \right]
     +n\Pr\left[\mathscr{D}_0^{l}\right], $$
where $\mathscr{D}_0^{l}$ is defined in Lemma \ref{lemma:rare}. For
$l \geq (\log x)^{3+\gamma} $, we
   obtain
     \begin{align*}
      \Pr\left[\rbp^{l} > l x\right]
       &= \Pr \left[ \sup_{j\geq 1}Z_{-j}>l x \right] \leq \Pr\left[\rbp_n^{l}>l x\right]
       +\Pr \left [\sup_{j > n}Z_{-j}>l x \right] \\
                          &\leq \Pr \left[\sup_{1\leq j \leq n}
                          \mp_j \left(1+ \underline{\mu}^{-1}
l^{-\frac{1}{3}} \right)^j >
                          x \right]+n\Pr\left[\mathscr{D}_0^{l}\right]+
                           \sum_{j>n}^{\infty}\Pr\left[Z_j^{l} >lx\right],
     \end{align*}
   which,  by setting $n=\lfloor h \log x
   \rfloor$ with $h$ being chosen as in Lemma  \ref{lemma:asymp_3} and
    applying Lemmas \ref{lemma:rare}, \ref{lemma:asymp_3},
yields,
\begin{align*}
 \Pr\left[\rbp^{l} > l x \right] &\leq \Pr \left[\sup_{1 \leq j \leq h \log x} \mp_j \left(1+\underline{\mu}^{-1}
l^{-\frac{1}{3}}\right)^j >
                            x \right]+o\left( \frac{1}{ x^{\alpha^{\ast} }} \right) \\
                          &\leq \Pr \left[\sup_{j \geq 1} \mp_j  >
                            x\left(1+\frac{1}
                          {\underline{\mu} (\log x)^{1+\gamma/3}} \right)^{-h \log x}\right]
                           +o\left( \frac{1}{ x^{\alpha^{\ast}}} \right).
     \end{align*}
 Finally, by using Theorem \ref{theorem:iidM} and observing that $\lim_{x\rightarrow \infty}
           \left(1+\frac{1}{\underline{\mu}(\log x)^{1+\gamma/3}} \right)^{-h \log x} =1$,
           we obtain
           $${\large \varlimsup_{\tiny \begin{array}{c}
  l\geq (\log x)^{3+\gamma} \\
  x\rightarrow \infty
\end{array}}} \Pr\left[\Lambda ^{l}/l>x \right]x^{\alpha^{\ast}}
         \leq \frac{1- \|G_{+}\| }{ \alpha^{\ast} \int_{0}^{\infty} xe^{\alpha^{\ast}x}G_{+}(dx)}. $$

Next, we prove the {lower bound}.
 Recall that $\mp_n^i \eqdef \prod_{j=i}^{n-1}\mu(J_j)$ and $\mp_i \eqdef
 \prod_{j=-1}^{-i}\mu(J_j)$. Then, for
   $l \geq (\log x)^{3+\gamma} $  and $n=\lfloor h \log x
   \rfloor$ where $h$ is chosen as in Lemma  \ref{lemma:asymp_3}, we have
 \begin{align}\label{eq:exactLower}
       \Pr\left[\rbp_n^l >lx\right] &\geq \Pr \left[ \sup_{0 \leq i \leq n-1}\mp_n^i \left(1-\underline{\mu}^{-1}
l^{-\frac{1}{3}} \right)^{n-i}> x \right] -
                    \Pr\left[\mathscr{E}_0^{l}\right]-\cdots-\Pr\left[\mathscr{E}_{n-1}^{l}\right] \nonumber\\
        &= \Pr \left[ \sup_{1 \leq i \leq n} \mp_i  \left(1-\underline{\mu}^{-1}
l^{-\frac{1}{3}} \right)^i> x \right] - n \Pr\left[\mathscr{E}_0^{l}\right] \nonumber\\
        &\geq \Pr \left[ \sup_{i \geq 1} \mp_i  \left(1-\underline{\mu}^{-1}
l^{-\frac{1}{3}} \right)^{h \log x}> x \right]
               -  \Pr \left[ \sup_{i > n} \mp_i  \left(1-\underline{\mu}^{-1}
l^{-\frac{1}{3}} \right)^i> x \right]
              - n \Pr\left[\mathscr{E}_0^{l}\right]\nonumber\\
      &\eqdef I_1-I_2-I_3,
    \end{align}
    where $\mathscr{E}_0^{l}$ is defined in Lemma \ref{lemma:rare}.
   By Lemma \ref{lemma:lowerI_2},  we obtain,
    \begin{align}\label{eq:exactLower2}
     I_2& \leq \Pr \left[
\sup_{i > h\log x} \mp_i > x \right]=o\left( \frac{1}{ x^{
\alpha^{\ast}} } \right),
\end{align}
      and by Lemma \ref{lemma:rare},
\begin{equation}\label{eq:exactLower3}
 I_3=o\left( \frac{1}{ x^{ \alpha^{\ast}} } \right).
\end{equation}
Thus, combining (\ref{eq:exactLower}), (\ref{eq:exactLower2}) and
(\ref{eq:exactLower3}), we obtain
 $$
    \Pr\left[\rbp_n^l >lx\right] \geq \Pr \left[\sup_{j \geq 1} \mp_j \left(1-\frac{1}
                          {\underline{\mu}\, (\log x)^{1+\gamma/3}} \right)^{h \log x} >
                            x\right]-o\left( \frac{1}{ x^{\alpha^{\ast} }} \right),
$$
which, by using the same argument as in the proof of the upper
bound, yields
 $${\large \varliminf_{\tiny \begin{array}{c}
  l\geq (\log x)^{3+\gamma} \\
  x\rightarrow \infty
\end{array}}} \Pr \left[\Lambda ^{l}/l>x \right]x^{\alpha^{\ast}}
         \geq \frac{1- \|G_{+}\| }{ \alpha^{\ast} \int_{0}^{\infty} xe^{\alpha^{\ast}x}G_{+}(dx)}. $$
\end{proof}

\subsection{Proofs of Theorems \ref{theorem:stopMul},  \ref{theorem:stopIID} and \ref{theorem:rsbp}}\label{ss:stopMulProof}
\begin{proof}[Proof of Theorem \ref{theorem:stopMul}]
First, we prove the \emph{upper bound}.  For a fixed $\alpha$ that
is in the neighborhood of $\alpha^{\ast}$ and $0<\epsilon<\lambda$,
there exists $n_{\epsilon}$ such that
     $\expect \left[ \left(\mp_n^0\right)^{\alpha} \right] < e^{(\Psi(\alpha)+\epsilon)n}$ and
     $e^{-(\lambda-\epsilon) n}>\Pr[N\geq n]> e^{-(\lambda+\epsilon) n}$
     for
 all
   $n \geq n_{\epsilon}$.  Since $\Psi(\alpha^{\ast})=\lambda$ and $\Psi'(\alpha^{\ast})>0$, we can choose $\delta, \epsilon>0$ small enough such that
  $\Psi(\alpha^{\ast}-\delta)-\lambda+2\epsilon=-\xi <0$.  Thus, noting that $N$ is independent of
   $\Pi_n$, we obtain
\begin{align}
\Pr\left[\Pi_N^0>x\right]&=\sum_{n=1}^{\infty}\Pr[N=n] \Pr\left[\Pi_n^0>x \right] \nonumber\\
&\leq \sum_{n=1}^{n_{\epsilon}}\Pr[N=n] \Pr\left[\Pi_n^0>x \right]+
         \sum_{n=n_{\epsilon}}^{\infty}\Pr[N \geq n] \Pr\left[\Pi_n^0>x \right] \nonumber\\
&\leq \sum_{n=1}^{n_{\epsilon}}\Pr[N=n]
\frac{\expect\left[\left(\Pi_n^0\right)^{\alpha^{\ast}}\right]}{x^{\alpha^{\ast}}}
     +\sum_{n=n_{\epsilon}}^{\infty} e^{-(\lambda-\epsilon) n}
     \frac{\expect[(\Pi_n^0 )^{\alpha^{\ast}-\delta}]}{x^{\alpha^{\ast}-\delta}}
     \nonumber\\
&\leq O\left(\frac{1}{x^{\alpha^{\ast}}}\right)
  +\frac{1}{x^{\alpha^{\ast}-\delta}} \sum_{n=n_{\epsilon}}^{\infty}
    e^{-\xi n}, \nonumber
\end{align}
which implies
$$
\varlimsup_{x\rightarrow \infty}\frac{\log
\Pr\left[\Pi_{N}^0>x\right]}{\log x} = -\alpha^{\ast}+\delta.
$$
Passing $\delta \to 0$ in the preceding equality completes the proof
of the upper bound.

Next, we prove the \emph{lower bound} by using the standard
exponential change of measure argument.
 For $0<3
\epsilon<\lambda$, $\delta>2\epsilon/(\lambda-3\epsilon)$ and $\log
x>n_{\epsilon}$, recalling that $e^{-(\lambda-\epsilon) n}>\Pr[N\geq
n]> e^{-(\lambda+\epsilon) n}$,  we obtain, for large $x$,
\begin{align*}
   \Pr\left[\frac{(1+\delta) \log x}{\Psi'(\alpha^{\ast})}\leq N
   \leq \frac{(1+2\delta) \log x}{\Psi'(\alpha^{\ast})}\right]
   &\geq e^{-\frac{(\lambda+\epsilon) (1+\delta) \log x}{\Psi'(\alpha^{\ast})}}-
        e^{-\frac{(\lambda-\epsilon)(1+2\delta) \log
   x}{\Psi'(\alpha^{\ast})}}\\
   &\geq (1-\epsilon)e^{-\frac{(\lambda+\epsilon) (1+\delta) \log x}{\Psi'(\alpha^{\ast})}}
\end{align*}
since $(\lambda+\epsilon) (1+\delta)<(\lambda-\epsilon)(1+2\delta)$,
and this implies that there exists $\delta \leq \zeta \leq 2\delta$
such that $n_x = \lceil (1+\zeta)(\log x)/\Psi'(\alpha^{\ast})
\rceil$ satisfies
\begin{equation}\label{eq:lowerN}
   \Pr[N=n_x]\geq \frac{(1-\epsilon)\Psi'(\alpha^{\ast}) e^{-(\lambda+\epsilon)(1+\delta) \log x /\Psi'(\alpha^{\ast})}}{\delta
   \log x}.
\end{equation}
Therefore, using (\ref{eq:lowerN}) and denoting $\log J_i$ by $X_i$,
we obtain
\begin{align}\label{eq:lowerN1}
 \Pr[\Pi_N^0>x ] &\geq \Pr[N=n_x] \Pr\left[\sum_{i=1}^{n_x}\log J_i> \log
 x
 \right] \nonumber\\
 &\geq \frac{(1-\epsilon)\Psi'(\alpha^{\ast})e^{-(\lambda+\epsilon)(1+\delta) \log x/\Psi'(\alpha^{\ast})}}{\delta
   \log
   x}\Pr\left[\sum_{i=1}^{n_x}X_i>\frac{\Psi'(\alpha^{\ast})}{1+\delta}n_x
 \right].
\end{align}

Next, we perform an exponential change of measure for the
probability on the right-hand side of (\ref{eq:lowerN1}). Let
$\Pr^{\ast}_n$ be the probability measure on $\real^n$ defined by
the probability measure $\Pr$ of the stationary and ergodic process
$\{X_i\}_{i\geq 1}$
$$
\Pr^{\ast}_n(dx_1,\cdots,dx_n)=e^{
\alpha^{\ast}S_n-\Psi_n(\alpha^{\ast})}\Pr(dx_1,\cdots,dx_n),
$$
where $S_n=\sum_{i=1}^{n}X_i$ and $\Psi_n(\alpha)\eqdef \log
\expect[e^{\alpha S_n}]$ satisfying $n^{-1}\Psi_n(\alpha) \to
\Psi(\alpha)$ in the neighborhood of $\alpha^{\ast}$.
Thus,
\begin{align}\label{eq:lowerN2}
\Pr\left[\sum_{i=1}^{n}X_i>\frac{\Psi'(\alpha^{\ast})}{1+\delta}n
 \right]&= \expect^{\ast}_n\left[ e^{-\alpha^{\ast}S_n+\Psi_n(\alpha^{\ast})}
      \ind \left(S_n>\frac{\Psi'(\alpha^{\ast})}{1+\delta}n \right) \right] \nonumber\\
      &\geq \expect^{\ast}_n\left[ e^{-\alpha^{\ast}S_n+\Psi_n(\alpha^{\ast})}
      \ind \left( \Big | \frac{S_n}{n}-\Psi'(\alpha^{\ast}) \Big | <
            \frac{\Psi'(\alpha^{\ast}) \delta}{1+\delta} \right) \right] \nonumber\\
           &\geq e^{-\alpha^{\ast}\frac{(1+2\delta)\Psi'(\alpha^{\ast})}{1+\delta}n+\Psi_n(\alpha^{\ast})}
      \Pr^{\ast}_n\left[ \Big | \frac{S_n}{n}-\Psi'(\alpha^{\ast}) \Big | <
            \frac{\Psi'(\alpha^{\ast}) \delta}{1+\delta} \right].
\end{align}
Then, by Claim $1$ on page $17$ of \cite{Buck90}, we know that
$$
 \Pr^{\ast}_n\left[ \Big | \frac{S_n}{n}-\Psi'(\alpha^{\ast}) \Big | <
            \frac{\Psi'(\alpha^{\ast}) \delta}{1+\delta} \right] \to 1\; \text{as $x \to \infty$},$$
            which, using (\ref{eq:lowerN1}) and
           setting $n=n_x$ in (\ref{eq:lowerN2}), yields
 $$
  \varliminf_{x \to \infty}\frac{\log \Pr\left[\Pi_N^0>x\right]}{\log
  x}\geq
  -\frac{(\lambda+\epsilon)(1+\delta)}{\Psi'(\alpha^{\ast})}- \frac{\alpha^{\ast}(1+2\delta)^2}{1+\delta}+
  \frac{(1+\delta)\Psi(\alpha^{\ast})}{\Psi'(\alpha^{\ast})}.
 $$
Finally, by passing $\epsilon, \delta \to 0$ in the preceding
equality and noting $\Psi(\alpha^{\ast})=\lambda$, we prove the
lower bound. \end{proof}

\begin{proof}[Proof of Theorem \ref{theorem:stopIID}]
 We give a constructive proof based on the connection (duality)
  between the $M/GI/1$ queue and the geometrically stopped multiplicative
 process.

Consider a $M/GI/1$ queue with the service distribution $\Pr[S\geq
t] = \bar{G}(t), t\geq 0$ and Poisson arrivals of rate $\lambda =
\rho/\expect\left[ S\right], \expect[S]<\infty$. Then, by the
\rm{Pollaczeck-Khinchine} formula (see, e.g., Theorem 5.7 on p.~237
of \cite{As87}), the stationary workload $Q$ of this $M/GI/1$ queue
is equal in distribution to $\sum_{i=1}^{N} H_i$, where $N,
\{H_i\}_{i\geq 1}$ are independent with $\Pr[N>n]=\rho^n, n\geq 0$
and
 \begin{equation*}
    \Pr[{H}_i\leq x] =
    \frac{\int_{0}^{x}\Pr[S\geq s]ds}{\expect [S]}=\frac{\int_{0}^{x}\bar{G}(s)ds}{\int_{0}^{\infty}\bar{G}(s)ds}
    =\Pr[\log J_i \leq x], \;\; x\geq 0,
 \end{equation*}
where the last equality follows from the assumption. Now, using the
preceding observation we show that there exists a RMP such that
$M=e^Q$ satisfies
 \begin{equation}\label{eq:stopIID}
  \Pr[M>x]=\Pr[Q>\log x] = \Pr \left[\sum_{i=1}^{N} H_i > \log x \right]
 =\Pr\left[\sum_{i=1}^{N} \log J_i > \log x \right]=\Pr\left[\Pi_{N}^0>
 x\right],
 \end{equation}
 which proves the first claim of the theorem.

  Next, using the
 additional assumptions of the theorem and applying Cram\'er-Lundberg theory for the $M/GI/1$ queue (e.g.,
see Theorem 5.2 in Chapter XIII of \cite{As87}), we obtain
 \begin{equation*}
   \lim_{x \rightarrow \infty} \Pr[M>
   x]x^{\alpha^{\ast}}=
 \frac{(1-\rho)\int_{0}^{\infty}\bar{G}(y)dy}{\alpha^{\ast}\rho \int_{0}^{\infty}ye^{\alpha^{\ast}y}\bar{G}(y)dy },
 \end{equation*}
which, by (\ref{eq:stopIID}), completes the proof.
\end{proof}

\begin{proof}[Proof of Theorem \ref{theorem:rsbp}]
The second equality is implied by Theorem~\ref{theorem:stopMul}, and
we only need to prove the first one. We begin with proving the
\emph{upper bound}. Recalling the definition of
$\mathcal{B}_n^{l,\epsilon}$ in Lemma \ref{lemma:B} and,  for $n\geq
1$, $0<\epsilon, \xi<1$, choosing $x^{\xi}>z_0>Z_0$, we obtain
\begin{align}
\Pr\left[Z_n^{Z_0}>x\right] & \leq \Pr\left[Z_n^{\lfloor x^{\xi}
\rfloor}>x\right]\nonumber\\
& \leq \Pr\left[Z_n^{\lfloor x^{\xi} \rfloor}>x,
\bigcap_{i=0}^{n-1}\left(\mathcal{B}_i^{\lfloor x^{\xi} \rfloor,\,
\epsilon}\right)^{C} \right] +
          \Pr\left[\bigcup_{i=0}^{n-1}\mathcal{B}_i^{\lfloor x^{\xi}
\rfloor ,\, \epsilon}\right] \nonumber \\
 &\leq \Pr\left[\Pi_n^0(1+\epsilon)^n>x^{1-\xi} \right] + n
\Pr\left[\mathcal{B}_0^{\lfloor x^{\xi}
\rfloor,\,\epsilon}\right],\nonumber
\end{align}
which, by the independence of $N$ and $\{B_n^i(j), J_n\}$, implies
   \begin{align}\label{eq:rsbp1}
   \Pr[Z_N>x]&\leq
\Pr\left[\Pi_N^0(1+\epsilon)^N>x^{1-\xi} \right] +
\expect[N]\Pr\left[\mathcal{B}_0^{\lfloor x^{\xi}
\rfloor,\,\epsilon}\right].
\end{align}

Next, define a new process $\{ \Pi_n^{\epsilon} = \Pi_n^0
(1+\epsilon)^n \}$. It is easy to see that, for
     $\epsilon$ small enough, the sequence $\{ \Pi_n^{\epsilon} \}$
     satisfies
      $ n^{-1}\log \expect \left[ (\Pi_n^{\epsilon})^{\alpha} \right]\rightarrow
      \Psi(\alpha)+\alpha \log(1+\epsilon)$. Therefore,
     by Theorem \ref{theorem:stopMul}, we obtain
  \begin{equation}\label{eq:rsbp2}
   \lim_{x\rightarrow \infty}\frac{\log \Pr\left[\Pi_N^0(1+\epsilon)^N>x^{1-\xi}\right]}{\log x}
= -(1-\xi)\alpha^{\ast}_{\epsilon},
\end{equation}
where $\alpha^{\ast}_{\epsilon}$ satisfies
$\Psi(\alpha^{\ast}_{\epsilon})+\alpha^{\ast}_{\epsilon}
\log(1+\epsilon)=0$. Combining (\ref{eq:rsbp1}), (\ref{eq:rsbp2})
and Lemma \ref{lemma:I_2}, we obtain
  \begin{equation*} \varlimsup_{x\rightarrow \infty}\frac{\log \Pr[Z_{N}>x ]}{\log x}
\leq -(1-\xi)\alpha^{\ast}_{\epsilon},
 \end{equation*}
 which, by passing $\epsilon, \xi \to 0$, completes the proof of the
 upper bound.

Now, we prove the \emph{lower bound}. Let $\{Z^1_{n,j}\}$ be i.i.d.
copies of $\{Z_n^1\}$ given the common modulating process $\{J_n\}$.
Then, noting that $Z^y_{n} \eqd \sum_{j=1}^{y}Z^1_{n,j}$ for integer
$y$ and using the union bound, we derive,  for $0<\xi<1, n\geq 0$,
\begin{align}
\Pr[Z_n>x] & \geq \frac{\lfloor x^{\xi} \rfloor }{ x^{\xi}
}\Pr[Z_n^1>x]\geq  \frac{1}{x^{\xi} }\Pr\left[Z_n^{\lfloor x^{\xi}
\rfloor}>x\lfloor x^{ \xi } \rfloor \right].\nonumber
\end{align}
Hence,  recalling the definition of $\mathcal{C}_n^{l,\epsilon}$ in
Lemma \ref{lemma:lowerI_3},  we  obtain
  \begin{align}
\Pr[Z_n>x] & \geq \frac{1}{x^{\xi}}\Pr\left[Z_n^{\lfloor x^{\xi}
\rfloor}>x\lfloor x^{ \xi } \rfloor,
\bigcap_{i=0}^{n-1}\left(\mathcal{C}_i^{\lfloor x^{\xi}
\rfloor,\,\xi}\right)^{C}
\right] \nonumber \\
&\geq \frac{1}{x^{\xi}} \left(\Pr\left[\Pi_n^0(1-\xi)^n>x \right] -
n \Pr\left[\mathcal{C}_0^{\lfloor x^{\xi} \rfloor,\,\xi}\right]
\right),\nonumber
\end{align}
which, by the independence of $N$ and $\{B_n^i(j), J_n\}$, yields
\begin{align}
\Pr[Z_N>x]&\geq  \frac{1}{x^{\xi} }\left(
\Pr\left[\Pi_N^0(1-\xi)^N>x \right] - \expect[N]
\Pr\left[\mathcal{C}_0^{\lfloor x^{\xi} \rfloor,\xi}\right] \right).
\nonumber
\end{align}
Then, by using the same approach as in the proof of the upper bound
and
 Lemma \ref{lemma:lowerI_3}, we can easily show that
  \begin{equation}
   \varliminf_{x\rightarrow \infty}\frac{\log \Pr[Z_{N}>x ]}{\log x}
\geq -\alpha^{\ast}. \nonumber
 \end{equation}

 Finally, by combining the upper bound and the lower bound, we finish the
 proof.
\end{proof}

\subsection{Proof of Theorem \ref{theorem:rdbranch}}\label{ss:absorbingProof}
\begin{proof}
 We begin with the \emph{upper bound}. Notice that when the system reaches
 stationarity,  $N_n$ follows the Poisson distribution, and therefore,  there exists $H>0$ such that
 \begin{equation}\label{eq:bpapN}
   \Pr[N_0>\lfloor H\log x \rfloor] = o\left( \frac{1}{x^{\alpha^{\ast}}} \right).
 \end{equation}
  Denoting by $\{\Lambda_{i}\}_{i\geq
 1}$ the i.i.d. copies of the random variable
 $\Lambda $ defined in Lemma \ref{lemma:zl},  we obtain
 \begin{align}
  \Pr[Z_s>x] &\leq \Pr\left[ \sum_{i=1}^{N_0}
  \Lambda_{i} >x \right] \nonumber\\
  &\leq \Pr\left[ \sum_{i=1}^{\lfloor H\log x \rfloor} \Lambda_{i} >x  \right]+
  \Pr[N_0 > \lfloor H \log x \rfloor ]\nonumber \\
  &\leq  H \log x \, \Pr\left[  \Lambda > \frac{x}{H\log x}  \right]+
  \Pr[N_0 > \lfloor H \log x  \rfloor  ], \nonumber
 \end{align}
 which, in conjunction with Theorem \ref{theorem:rbp} and equation
 (\ref{eq:bpapN}), yields
  \begin{equation}\label{eq:bpabUpper}
\varlimsup_{x \rightarrow \infty}\frac{\log \Pr[Z_s>x]}{\log x} \leq
-\alpha^{\ast}.
 \end{equation}

 Next, we proceed with the \emph{lower bound}. Construct a new process
 that has the same arrivals $\{A_n\}$ as described before
  but only allows at most one object to exist
 in the system.   The
construction goes as follows: all the new arrivals will be dropped
if there is an object present in the system;
  similarly, when newly generated objects arrive to the empty system, only one object
 will be accepted while others will be dropped;  the object, if any,  evolves according to an
i.i.d. copy of the modulated branching process $Z_P$. Denote the
total size of the object in the new system at time $n$ by
$\underline{Z}_n$, and observe that $\underline{Z}_n$ forms a
renewal process. Then, taking out all the empty (idle) periods of
the new system and concatenating the remaining periods sequentially
yields a process equal in distribution to a reflected modulated
branching process $\{\Lambda_n\}$, as defined in
(\ref{eq:rbranching}). Therefore, when the new system is in
stationarity, we obtain, by the independence of $\{A_n\}$ and $Z_P$,
\begin{align}
 \Pr[Z_s>x]&\geq \Pr[\underline{Z}_0>x]
 = \Pr[\underline{Z}_0>0]\Pr[\Lambda>x], \nonumber
\end{align}
which, by Theorem \ref{theorem:rbp}, yields
\begin{equation}\label{eq:bpabLower}
\varliminf_{x \rightarrow \infty}\frac{\log \Pr[Z_s>x]}{\log x} \geq
-\alpha^{\ast}.
 \end{equation}

 Finally, combining (\ref{eq:bpabUpper}) and (\ref{eq:bpabLower})
 finishes the proof.
 \end{proof}

\section*{Acknowledgements}
We are thankful to the anonymous reviewers for their helpful
comments and, in particular, we are grateful to a reviewer who
pointed to us reference \cite{Goldie91} that led to Theorem
\ref{theorem:ExactMBP}.
 \small
\bibliography{JelenkovicTanMBP08}

\end{document}